\documentclass[12pt]{amsart}
\usepackage{amssymb,amscd}
\usepackage{pstricks}

\newtheorem{theorem}{Theorem}[section]
\newtheorem{lemma}[theorem]{Lemma}
\newtheorem{proposition}[theorem]{Proposition}
\newtheorem{corollary}[theorem]{Corollary}
\newtheorem{definition}[theorem]{Definition}
\newtheorem{conj}[theorem]{Conjecture}

\newtheorem{assumption}[theorem]{Assumption}
\theoremstyle{definition}
\newtheorem{remark}[theorem]{Remark}

\numberwithin{equation}{section}

\newcommand{\Z}{\mathbb{Z}}
\newcommand{\R}{\mathbb{R}}
\newcommand{\C}{\mathbb{C}}
\newcommand{\CP}{\mathbb{CP}}

\renewcommand{\Re}{\mathrm{Re}\,}
\renewcommand{\Im}{\mathrm{Im}\,}
\newcommand{\hol}{\mathrm{hol}}
\newcommand{\LLog}{\text{\rm\L{}og}}
\newcommand{\m}{\mathfrak{m}}

\title[Mirror symmetry and T-duality]{Mirror symmetry and T-duality in the complement of an
anticanonical divisor}
\author{Denis Auroux}
\thanks{Partially supported by NSF grant DMS-0600148 and an A.P. Sloan
research fellowship.}
\address{Department of Mathematics, M.I.T., Cambridge MA 02139, USA}
\email{auroux@math.mit.edu}

\begin{document}
\begin{abstract}
We study the geometry of complexified moduli spaces of special
Lagrangian submanifolds in the complement of an anticanonical divisor
in a compact K\"ahler manifold. In particular, we explore the connections
between T-duality and mirror symmetry in concrete examples, and show how
quantum corrections arise in this context.
\end{abstract}

\maketitle

\section{Introduction}

The Strominger-Yau-Zaslow conjecture \cite{SYZ} asserts that the mirror
of a Calabi-Yau manifold can be constructed by dualizing a fibration by
special Lagrangian tori. This conjecture has been studied extensively,
and the works of Fukaya, Kontsevich and Soibelman, Gross and Siebert, and
many others paint a very rich and subtle picture of mirror symmetry as a
T-duality modified by ``quantum corrections'' \cite{GS,KS1,KS2}.

On the other hand, mirror symmetry has been extended to the non Calabi-Yau
setting, and in particular to Fano manifolds, by considering Landau-Ginzburg
models, i.e.\ noncompact manifolds equipped with a complex-valued function
called {\it superpotential} \cite{HV}. Our goal is to understand the
connection between mirror symmetry and T-duality in this setting. 

For a toric Fano manifold, the moment map provides a fibration by Lagrangian
tori, and in this context the mirror construction can be understood as a
T-duality, as evidenced e.g.\ by Abouzaid's work \cite{Ab1,Ab2}. Evidence
in the non-toric case is much scarcer, in spite of Hori and Vafa's
derivation of the mirror for Fano complete intersections in toric varieties
\cite{HV}. The best understood case so far is that of Del Pezzo surfaces
\cite{AKO2}; however, in that example the construction of the mirror is
motivated by entirely ad hoc considerations. As an attempt to understand
the geometry of mirror symmetry beyond the Calabi-Yau setting, we start
by formulating the following naive conjecture:

\begin{conj}\label{conj:main}
Let $(X,\omega,J)$ be a compact K\"ahler manifold, let $D$ be an
anticanonical divisor in $X$, and let $\Omega$ be a holomorphic volume
form defined over $X\setminus D$. Then a mirror manifold $M$ can be
constructed as a moduli space of special Lagrangian tori in $X\setminus D$
equipped with flat\/ $U(1)$ connections over them, with a
superpotential\/ $W:M\to\C$ given by Fukaya-Oh-Ohta-Ono's $m_0$ obstruction
to Floer homology. Moreover, the fiber of this Landau-Ginzburg model is
mirror to $D$.
\end{conj}

The main goal of this paper is to investigate the picture
suggested by this conjecture. Conjecture \ref{conj:main}
cannot hold as stated, for several reasons. One is that
in general the special Lagrangian torus fibration on $X\setminus D$ is
expected to have singular fibers, which requires suitable corrections to
the geometry of $M$. Moreover, the superpotential
constructed in this manner is not well-defined, since wall-crossing
phenomena make $m_0$ multivalued. In particular it is not clear how to
define the fiber of $W$. These various issues are related to quantum
corrections arising from holomorphic discs of Maslov index 0; while we do
not attempt a rigorous systematic treatment, general considerations 
(see \S \ref{ss:wall1}--\ref{ss:wall2}) and calculations
on a specific example (see Section \ref{s:nontoric}) suggest that the
story will be very similar to the Calabi-Yau case \cite{GS,KS2}.
Another issue is the incompleteness of $M$; according to Hori and
Vafa \cite{HV}, this is an indication that the mirror symmetry
construction needs to be formulated in a certain renormalization limit
(see \S \ref{ss:compare}). The modifications of Conjecture \ref{conj:main}
suggested by these observations are summarized in Conjectures
\ref{conj:quantumcorr} and \ref{conj:renormalize} respectively.

The rest of this paper is organized as follows.
In Section \ref{s:defs} we
study the moduli space of special Lagrangians and its geometry. In Section
\ref{s:W} we discuss the $m_0$ obstruction in Floer theory and the
superpotential. Then Section \ref{s:toric} is devoted to the toric case
(in which the superpotential was already investigated by Cho and Oh \cite{ChoOh}),
and Section \ref{s:nontoric} discusses in detail the example of $\CP^2$
with a non-toric holomorphic volume form. Finally, Section \ref{s:m0c1}
explores the relation between the critical values of $W$ and
the quantum cohomology of $X$, and Section \ref{s:fiber} discusses the
connection to mirror symmetry for the Calabi-Yau hypersurface $D\subset X$.

Finally, a word of warning is in order: in the interest of readability and
conciseness, many of the statements made in this paper are not entirely
rigorous; in particular, weighted counts of holomorphic discs are always
assumed to be convergent, and issues related to the lack of regularity of
multiply covered Maslov index 0 discs are mostly ignored. Since the main
goal of this paper is simply to evidence specific phenomena and illustrate
them by examples, we feel that this approach is not unreasonable, and ask
the detail-oriented reader for forgiveness.

\subsection*{Acknowledgements}
I am heavily indebted to Mohammed Abouzaid, Paul Seidel and Ludmil Katzarkov
for numerous discussions which played a crucial role in the
genesis of this paper. I would also like to thank Leonid Polterovich and
Felix Schlenk for their explanations concerning the Chekanov torus, as
well as Anton Kapustin, Dima Orlov and Ivan Smith for helpful discussions.
This work was partially supported by an NSF grant (DMS-0600148) and an
A.P. Sloan research fellowship.

\section{The complexified moduli space of special Lagrangians}\label{s:defs}
\subsection{Special Lagrangians}
Let $(X,\omega,J)$ be a smooth compact K\"ahler manifold of complex
dimension $n$, and let $\sigma\in H^0(X,K_X^{-1})$ be a nontrivial
holomorphic section of the anticanonical bundle, vanishing on a divisor
$D$. Then the complement $X\setminus D$ carries a nonvanishing holomorphic 
$n$-form $\Omega=\sigma^{-1}$. By analogy with the Calabi-Yau situation,
for a given $\phi\in \R$ we make the following definition:

\begin{definition}
A Lagrangian submanifold $L\subset X\setminus D$ is {\em special Lagrangian}
with phase $\phi$ if\/ $\Im(e^{-i\phi}\Omega)_{|L}=0$.
\end{definition}

Multiplying $\Omega$ by $e^{-i\phi}$ if necessary, in the rest of this
section we will consider the case $\phi=0$. In the Calabi-Yau case,
McLean has shown that infinitesimal deformations of special Lagrangian
submanifolds correspond to harmonic 1-forms, and that these deformations
are unobstructed \cite{mclean}. (See also \cite{hitchin} and \cite{salur}
for additional context).

In our case, the restriction to $L$ of
$\Re(\Omega)$ is a non-degenerate volume form (which we assume
to be compatible with the orientation of $L$), but it differs from
the volume form $vol_g$ induced by the K\"ahler metric $g$. Namely, there
exists a function $\psi\in C^\infty(L,\R_+)$ such that
$\Re(\Omega)_{|L}=\psi\,vol_g$.

\begin{definition}
A one-form $\alpha\in \Omega^1(L,\R)$ is {\em $\psi$-harmonic} if
$d\alpha=0$ and $d^*(\psi\alpha)=0$. We denote by $\mathcal{H}^1_\psi(L)$
the space of $\psi$-harmonic one-forms.
\end{definition}

\begin{lemma}
Each cohomology class contains a unique $\psi$-harmonic representative.
\end{lemma}

\proof
If $\alpha=df$ is exact and $\psi$-harmonic, then $\psi^{-1}d^*(\psi\,df)=
\Delta f-\psi^{-1}\langle d\psi,df\rangle=0$. Since the maximum principle
holds for solutions of this equation, $f$ must be constant. So
every cohomology class contains at most one $\psi$-harmonic representative.

To prove existence, we consider the elliptic operator
$D:\Omega^\mathrm{odd}(L,\R)\to \Omega^\mathrm{even}(L,\R)$
defined by $D(\alpha_1,\alpha_3,\dots)=(\psi^{-1}d^*(\psi\alpha_1),
d\alpha_1+d^*\alpha_3,\dots)$.  Clearly the kernel of
$D$ is spanned by $\psi$-harmonic 1-forms and by harmonic forms of odd
degree $\ge 3$, while its cokernel contains all harmonic forms of even
degree $\ge 2$ and the function $\psi$. However $D$ differs from $d+d^*$
by an order 0 operator, so its index is
$\mathrm{ind}(D)=\mathrm{ind}(d+d^*)=-\chi(L)$. It follows that $\dim
\mathcal{H}^1_\psi(L)=\dim H^1(L,\R)$.
\endproof

\begin{remark}
Rescaling the metric by a factor of $\lambda^2$ modifies the Hodge $*$
operator on 1-forms by a factor of $\lambda^{n-2}$. Therefore,
if $n\neq 2$, then a 1-form is $\psi$-harmonic if and only if it is
harmonic for the rescaled metric $\tilde{g}=\psi^{2/(n-2)}g$.
\end{remark}

\begin{proposition} \label{prop:deform}
Infinitesimal special Lagrangian deformations of $L$ are in one to one
correspondence with $\psi$-harmonic 1-forms on $L$.
More precisely, a section of the normal bundle $v\in C^\infty(NL)$ determines
a 1-form $\alpha=-\iota_v\omega\in \Omega^1(L,\R)$ and an \hbox{$(n-1)$}-form
$\beta=\iota_v\Im\Omega\in \Omega^{n-1}(L,\R)$. These satisfy
\hbox{$\beta=\psi\,*_g\alpha$}, and the
deformation is special Lagrangian if and only if $\alpha$ and $\beta$ are
both closed.
Moreover, the deformations are unobstructed.
\end{proposition}

\proof
For special Lagrangian $L$, we have linear isomorphisms
$NL\simeq T^*L\simeq \wedge^{n-1}T^*L$ given
by the maps $v\mapsto -\iota_v\omega$ and 
$v\mapsto \iota_v\Im\Omega$.
More precisely, given a point $p\in L$, by complexifying a $g$-orthonormal
basis of $T_pL$ we obtain a local frame $(\partial_{x_j},\partial_{y_j})$
in which $\omega,J,g$ are standard at $p$, and
$T_pL=\mathrm{span}(\partial_{x_1},\dots,\partial_{x_n})$. In terms of
the dual basis $dz_j=dx_j+idy_j$, at the point $p$ we have $\Omega=\psi\,
dz_1\wedge\dots\wedge dz_n$. Hence, given $v=\sum c_j\partial_{y_j}
\in N_pL$, we have $-\iota_v\omega=\sum c_j\,dx_j$ and $$\textstyle \iota_v\Im\Omega=
\psi\, \sum\limits_j c_j (-1)^{j-1}dx_1\wedge\dots\wedge \widehat{dx_j}\wedge\dots\wedge
dx_n=\psi\,*_g(-\iota_v\omega).$$
Consider a section of the normal bundle $v\in C^\infty(NL)$, and use an
arbitrary metric to construct a family of submanifolds $L_t=j_t(L)$, where
$j_t(p)=\exp_p(tv(p))$. Since $\omega$ and $\Im\Omega$ are closed,
we have $$\frac{d}{dt}_{|t=0}(j_t^*\omega)=L_v\omega=
d(\iota_v\omega)\quad\mathrm{and}\quad
\frac{d}{dt}_{|t=0}(j_t^*\Im\Omega)=
L_v\Im\Omega=d(\iota_v\Im\Omega).$$
Therefore, the infinitesimal deformation $v$ preserves the special
Lagrangian condition $\omega_{|L}=\Im\Omega_{|L}=0$ if and only if
the forms $\alpha=-\iota_v\omega$ and $\beta=\iota_v\Im\Omega$ are closed.
Since $\beta=\psi\,*_g\alpha$, this is equivalent to the requirement
that $\alpha$ is $\psi$-harmonic.

Finally, unobstructedness is proved exactly as in the Calabi-Yau case,
by observing that the linear map $v\mapsto (L_v\omega,L_v\Im\Omega)$ 
from normal vector fields to exact 2-forms and exact $n$-forms is
surjective and invoking the implicit function theorem \cite{mclean}.
\endproof

This proposition allows us to consider (at least locally) the moduli space
of special Lagrangian deformations of $L$. This moduli space is a smooth
manifold, and carries two
natural integer affine structures, obtained by identifying the tangent space
to the moduli space with either $H^1(L,\R)$ or $H^{n-1}(L,\R)$ and
considering the integer cohomology lattices. 

\subsection{The geometry of the complexified moduli space}

In this section we study the geometry of the (complexified) moduli space
of special Lagrangian submanifolds. In the Calabi-Yau case, our
constructions essentially reduce to those in Hitchin's illuminating paper
\cite{hitchin}.

We now consider pairs $(L,\nabla)$ consisting of a special Lagrangian
submanifold $L\subset X\setminus D$ and a flat unitary connection $\nabla$
on the trivial complex line bundle over $L$, up to gauge equivalence.
(In the presence of a B-field we would instead require $\nabla$ to have
curvature $-iB$; here we do not consider B-fields). Allowing $L$ to vary
in a given $b_1(L)$-dimensional family $\mathcal{B}$ of special Lagrangian
submanifolds (a domain in the moduli space), we
denote by $M$ the space of equivalence classes of pairs $(L,\nabla)$.
Our first observation is that $M$ carries a natural integrable complex structure.

Indeed, recall that the gauge equivalence
class of the connection $\nabla$ is determined by its holonomy 
$\hol_\nabla\in \mathrm{Hom}(H_1(L), U(1))\simeq H^1(L,\R)/H^1(L,\Z)$.
We will choose a representative of the form $\nabla=d+iA$,
where $A$ is a $\psi$-harmonic 1-form on $L$.

Then the tangent space to $M$ at
a point $(L,\nabla)$ is the set of all pairs $(v,\alpha)\in C^\infty(NL)
\oplus \Omega^1(L,\R)$ such that
$v$ is an infinitesimal special Lagrangian deformation, and
$\alpha$ is a $\psi$-harmonic 1-form, viewed as an infinitesimal
deformation of the flat connection. The map $(v,\alpha)\mapsto -\iota_v\omega+
i\alpha$ identifies $T_{(L,\nabla)}M$ with the space
$\mathcal{H}^1_{\psi} (L)\otimes \C$ of complex-valued
$\psi$-harmonic 1-forms on $L$, which makes $M$ a complex manifold.
More explicitly, the complex structure on $M$ is as follows:

\begin{definition}
Given $(v,\alpha)\in T_{(L,\nabla)}M\subset C^\infty(NL)\oplus
\Omega^1(L,\R)$, we define
$J^\vee(v,\alpha)=(a,-\iota_v\omega)$, where $a$ is
the normal vector field such that $\iota_a\omega=\alpha$.
\end{definition}

The following observation will be useful in Section \ref{s:W}:

\begin{lemma}\label{l:holomcoords}
Let $A\in H_2(M,L;\Z)$ be a relative homology class with boundary $\partial
A\neq 0 \in H_1(L,\Z)$. Then the function
\begin{equation}\label{eq:coord}
z_A= \exp(-\textstyle\int_A\omega)\,\hol_\nabla(\partial A):M\to\C^*
\end{equation}
is holomorphic.
\end{lemma}

\proof The differential $d\log z_A$ is simply $(v,\alpha)\mapsto
\int_{\partial A}-\iota_v\omega+i\alpha$, which is $\C$-linear.\endproof

More precisely, the function $z_A$ is well-defined locally (as long as we can
keep track of the relative homology class $A$ under deformations of $L$),
but might be multivalued if the family of special Lagrangian deformations 
of $L$ has non-trivial monodromy.

If the map $j_*:H_1(L)\to H_1(X)$ induced by inclusion is trivial, then this
yields a set of (local) holomorphic coordinates $z_i=z_{A_i}$
on $M$, by considering a collection of relative homology classes $A_i$
such that $\partial A_i$ form a basis of $H_1(L)$. Otherwise, given a class
$c\in H_1(L)$ we can fix a representative $\gamma^0_c$ of the class
$j_*(c)\in H_1(X)$, and use the symplectic area of a 2-chain in $X$ with
boundary on $\gamma^0_c\cup L$, together with the holonomy of $\nabla$
along the part of the boundary contained in $L$, as a substitute for the
above construction.

Next, we equip $M$ with a symplectic form:
\begin{definition}\label{def:omegavee}
Given $(v_1,\alpha_1),(v_2,\alpha_2)\in T_{(L,\nabla)}M$, we define
$$\omega^\vee((v_1,\alpha_1),(v_2,\alpha_2))=\int_L \alpha_2\wedge
\iota_{v_1}\Im\Omega-\alpha_1\wedge \iota_{v_2}\Im\Omega.$$
\end{definition}

\begin{proposition}
$\omega^\vee$ is a K\"ahler form on $M$, compatible with $J^\vee$.
\end{proposition}

\proof
First we prove that $\omega^\vee$ is closed and non-degenerate by
exhibiting local coordinates on $M$ in which it is standard.
Let $\gamma_1,\dots,\gamma_r$ be a basis of $H_{n-1}(L,\Z)$ (modulo
torsion), and let $e^1,\dots,e^r$ be the Poincar\'e dual basis of
$H^1(L,\Z)$. Let $\gamma^1,\dots,\gamma^r$ and $e_1,\dots,e_r$ be the
dual bases of $H^{n-1}(L,\Z)$ and $H_1(L,\Z)$ (modulo torsion): then
$\langle e^i\cup \gamma^j,[L]\rangle = \langle \gamma^j,
\gamma_i\rangle = \delta_{ij}$. In particular,
for all $a\in H^1(L,\R)$ and $b\in H^{n-1}(L,\R)$
we have \begin{equation}\label{eq:poincaredual}\textstyle
\langle a\cup b,
[L]\rangle=\sum\limits_{i,j} \langle a,e_i\rangle \langle
b,\gamma_j\rangle \langle e^i\cup \gamma^j,[L]\rangle=
\sum\limits_i \langle a,e_i\rangle \langle
b,\gamma_i\rangle.\end{equation}
Fix representatives $\Gamma_i$ and $E_i$ of the homology classes
$\gamma_i$ and $e_i$, and 
consider a point $(L',\nabla')$ of $M$ near $(L,\nabla)$. $L'$ is the image
of a small deformation $j'$ of the inclusion map $j:L\to X$.
Consider an $n$-chain $C_i$ in $X\setminus D$ such that $\partial C_i=
j'(\Gamma_i)-j(\Gamma_i)$, and let $p_i=\int_{C_i}\Im\Omega$. Also, let $\theta_i$
be the integral over $E_i$ of the connection 1-form of $\nabla'$ in a fixed
trivialization. Then $p_1,\dots,p_r,\theta_1,\dots,\theta_r$ are local coordinates
on $M$ near $(L,\nabla)$, and their differentials are given by $dp_i(v,\alpha)=
\langle [\iota_v\Im \Omega],\gamma_i\rangle$ and $d\theta_i(v,\alpha)=
\langle [\alpha],e_i\rangle$. 
Using (\ref{eq:poincaredual}) we deduce that $\omega^\vee=\sum_{i=1}^r
dp_i\wedge d\theta_i$.

Next we observe that, by Proposition \ref{prop:deform}, 
$\omega^\vee((v_1,\alpha_1),(v_2,\alpha_2))$ can be rewritten as 
$$\int_L \alpha_1\wedge (\psi\,*\!\iota_{v_2}\omega)-\alpha_2\wedge (\psi\,*\!
\iota_{v_1}\omega)=\int_L \psi\,\bigl(\langle \alpha_1,\iota_{v_2}\omega\rangle_g
-\langle \iota_{v_1}\omega,\alpha_2\rangle_g\bigr)\,vol_g.$$
So the compatibility of $\omega^\vee$ with $J^\vee$ follows directly from
the observation that
$$\omega^\vee((v_1,\alpha_1),J^\vee(v_2,\alpha_2))=
\int_L \psi\,\bigl(\langle \alpha_1,\alpha_2\rangle_g+\langle
\iota_{v_1}\omega,\iota_{v_2}\omega\rangle_g\bigr)\,vol_g$$
is clearly a Riemannian metric on $M$.
\endproof

\begin{remark}
Consider the projection $\pi:M\to \mathcal{B}$ which forgets the
connection, i.e.\ the map $(L,\nabla)\mapsto L$. Then the
fibers of $\pi$ are Lagrangian with respect to $\omega^\vee$. 
\end{remark}

If $L$ is a torus, then $\dim M=\dim X=n$ and we can also equip $M$ with
a holomorphic volume form defined as follows:
\begin{definition}
Given $n$ vectors $(v_1,\alpha_1),\dots,(v_n,\alpha_n)\in T_{(L,\nabla)}M
\subset C^\infty(NL)\oplus \Omega^1(L,\R)$, we define
$$\Omega^\vee((v_1,\alpha_1),\dots,(v_n,\alpha_n))=\int_L
(-\iota_{v_1}\omega+i\alpha_1)\wedge\dots\wedge
(-\iota_{v_n}\omega+i\alpha_n).$$
\end{definition}
In terms of the local holomorphic coordinates $z_1,\dots,z_n$ on $M$
constructed from a basis of $H_1(L,\Z)$ using the discussion after
Lemma \ref{l:holomcoords}, this holomorphic volume form is simply
$d\log z_1\wedge\dots\wedge d\log z_n$.

In this situation, the fibers of $\pi:M\to\mathcal{B}$ are special
Lagrangian (with phase $n\pi/2$) with respect to $\omega^\vee$ and
$\Omega^\vee$.
If in addition we assume that $\psi$-harmonic 1-forms on $L$ have no 
zeroes (this is automatic in dimensions $n\le 2$ using the maximum
principle), then we recover the familiar picture:
in a neighborhood of $L$, $(X,J,\omega,\Omega)$ and $(M,J^\vee,
\omega^\vee,\Omega^\vee)$ carry dual fibrations by special Lagrangian
tori.

\section{Towards the superpotential}\label{s:W}

\subsection{Counting discs}

Thanks to the monumental work of Fukaya, Oh, Ohta and Ono \cite{FO3},
it is now well understood that the Floer complex of a Lagrangian
submanifold carries the structure of a {\em curved}\/ or {\em obstructed}
$A_\infty$-algebra. The key ingredient is the moduli space of
$J$-holomorphic discs with boundary in the given Lagrangian submanifold,
together with evaluation maps at boundary marked points. In our case
we will be mainly interested in (weighted) counts of holomorphic discs
of Maslov index 2 whose boundary passes through a given point of the
Lagrangian; in the Fukaya-Oh-Ohta-Ono formalism, this corresponds to the
degree 0 part of the obstruction term $\m_0$. In the toric case it is 
known that this quantity agrees with the superpotential
of the mirror Landau-Ginzburg model; see in particular
the work of Cho and Oh \cite{ChoOh}, and
\S \ref{s:toric} below. In fact, the material in this section
overlaps signiicantly with \cite{ChoOh}, and with \S 12.7 of \cite{FO3}.

As in \S \ref{s:defs}, we consider a smooth compact K\"ahler manifold
$(X,\omega,J)$ of complex dimension $n$, equipped with a holomorphic
$n$-form $\Omega$ defined over the complement of an anticanonical divisor
$D$. 

Recall that, given a Lagrangian submanifold $L$ and a nonzero relative homotopy
class $\beta\in \pi_2(X,L)$, the moduli space $\mathcal{M}(L,\beta)$ of
$J$-holomorphic discs with boundary on $L$ representing the class $\beta$
has virtual dimension $n-3+\mu(\beta)$, where $\mu(\beta)$ is the Maslov
index.

\begin{lemma}\label{l:maslov}
If $L\subset X\setminus D$ is special Lagrangian, then $\mu(\beta)$ is
equal to twice the algebraic intersection number $\beta\cdot[D]$.
\end{lemma}

\proof
Because the tangent space to $L$ is totally real, the choice of a volume
element on $L$ determines a nonvanishing section $\det(TL)$ of $K_X^{-1}=
\Lambda^n(TX,J)$ over $L$. Its square $\det(TL)^{\otimes 2}$ defines a section of
the circle bundle $S(K_X^{-2})$ associated to $K_X^{-2}$ over $L$,
independent of the chosen volume element.
The Maslov number $\mu(\beta)$ measures the obstruction of this section 
to extend over a disc $\Delta$ representing the class $\beta$
(see Example 2.9 in \cite{SeGraded}). 

Recall that $D$ is the divisor
associated to $\sigma=\Omega^{-1}\in H^0(X,K_X^{-1})$. Then
$\sigma^{\otimes 2}$ defines a section of $S(K_X^{-2})$
over $L\subset X\setminus D$, and since $L$ is special Lagrangian,
the sections $\sigma^{\otimes 2}$ and $\det(TL)^{\otimes 2}$ coincide
over $L$ (up to a constant phase factor $e^{-2i\phi}$).
Therefore, $\mu(\beta)$ measures precisely the obstruction for
$\sigma^{\otimes 2}$ to extend over $\Delta$, which is twice
the intersection number of $\Delta$ with $D$.
\endproof

In fact, as pointed out by M.\ Abouzaid, the same result holds if we
replace the special Lagrangian condition by the weaker requirement that
the Maslov class of $L$ vanishes in $X\setminus D$ (i.e., 
the phase function $\arg(\Omega_{|L})$ lifts to a real-valued function).

Using positivity of intersections, Lemma \ref{l:maslov} implies that
all holomorphic discs with boundary in $L$ have non-negative Maslov index.

We will now make various assumptions on $L$ in order to ensure that the count
of holomorphic discs that we want to consider is well-defined:

\begin{assumption}\leavevmode\label{as:good}
\begin{enumerate}
\item there are no non-constant holomorphic discs of Maslov index 0 in $(X,L)$;
\item holomorphic discs of Maslov index 2 in $(X,L)$ are regular;
\item there are no non-constant holomorphic spheres in $X$ with $c_1(TX)\cdot[S^2]\le 0$.
\end{enumerate}
\end{assumption}

Then, for every relative homotopy class $\beta\in \pi_2(X,L)$ such that
$\mu(\beta)=2$, the moduli space $\mathcal{M}(L,\beta)$
of holomorphic discs with boundary in $L$
representing the class $\beta$ is a smooth compact manifold of real 
dimension $n-1$: no bubbling or multiple covering phenomena can occur
since 2 is the minimal Maslov index. 

We also assume that $L$ is spin (recall that we are chiefly
interested in tori), and choose a spin structure. The choice is
not important, as the difference between two spin structures is an element
of $H^1(L,\Z/2)$ and can be compensated by twisting the connection $\nabla$
accordingly. Then $\mathcal{M}(L,\beta)$ is oriented, and 
the evaluation map at a boundary marked point gives us an $n$-cycle
in $L$, which is of the form $n_\beta(L)\,[L]$ for some integer
$n_\beta(L)\in\Z$. In simpler terms, $n_\beta(L)$ is the (algebraic) number of holomorphic discs in
the class $\beta$ whose boundary passes through a generic point $p\in L$.

Then, ignoring convergence issues, we can tentatively make the following
definition (see also \cite{ChoOh},
\S 12.7 in \cite{FO3}, and Section 5b in \cite{seidel}):

\begin{definition} \label{def:m0}
$\displaystyle m_0(L,\nabla)=\sum_{\beta,\ \mu(\beta)=2} n_\beta(L)\,
\exp(-\textstyle\int_\beta\omega)\,\hol_\nabla(\partial\beta).$
\end{definition}

If Assumption \ref{as:good} holds for all special Lagrangians in the
considered family $\mathcal{B}$, and if the sum converges, then we obtain
in this way a complex-valued function on $M$, which we call {\it superpotential}
and also denote by $W$ for consistency with the literature.
In this ideal situation, the integers $n_\beta(L)$ are locally constant,
and Lemma \ref{l:holomcoords} immediately implies:

\begin{corollary}
$W=m_0:M\to\C$ is a holomorphic function.
\end{corollary}

An important example is the case of toric fibers in a toric manifold,
discussed in Cho and Oh's work \cite{ChoOh} and in \S \ref{s:toric} below:
in this case, the superpotential $W$ agrees with Hori and Vafa's physical
derivation \cite{HV}.

\begin{remark} \label{rmk:formal}
The way in which we approach the superpotential here is a bit different
from that in \cite{FO3}. Fukaya, Oh, Ohta and Ono consider a single
Lagrangian submanifold $L$, and the function which to a 1-cocycle $a$
associates the degree zero part of $\m_0+\m_1(a)+\m_2(a,a)+\dots$. However,
each of these terms counts holomorphic discs of Maslov index 2 whose boundary
passes through a generic point of $L$, just with different weights.
It is not hard to convince oneself that the contribution to $\m_k(a,a,\dots)$
of a disc in a given
class $\beta$ is weighted by a factor $\frac{1}{k!}\langle a,\partial\beta
\rangle^k$ (the coefficient $\frac{1}{k!}$ comes from the requirement that
the $k$ input marked points must lie in the correct order on the boundary of
the disc). Thus, the series $\m_0+\m_1(a)+\m_2(a,a)+\dots$ counts
Maslov index 2 discs with weights $\exp(\int_{\partial \beta} a)$ (in
addition to the weighting by symplectic area). In this sense $a$ can
be thought of as a {\it non-unitary} holonomy (normally with values
in the positive part of the Novikov ring for convergence reasons; here
we assume convergence and work with complex numbers). Next, we observe that,
since the weighting by symplectic area and holonomy is encoded by the
complex parameter $z_\beta$ defined in (\ref{eq:coord}), varying the
holonomy in a non-unitary manner is equivalent to moving the Lagrangian
in such a way that the flux of the symplectic form equals the
real part of the connection form. More precisely, this equivalence between
a non-unitary connection on a fixed $L$ and a unitary connection on a
non-Hamiltonian deformation of $L$ only holds as long as the disc
counts $n_\beta$ remain constant; so in general the superpotential in
\cite{FO3} is the analytic continuation of the germ of our superpotential
at the considered point.
\end{remark}

\begin{remark} \label{rmk:negspheres}
Condition (3) in Assumption \ref{as:good} can be somewhat relaxed.
For example, one can allow the existence of nonconstant $J$-holomorphic
spheres of Chern number 0, as long as all simple (non multiply covered)
such spheres are regular, and the associated evaluation
maps are transverse to the evaluation maps at interior marked points of
$J$-holomorphic discs of Maslov index 2 in $(X,L)$. Then the union of all
holomorphic spheres with Chern number zero is a subset $\mathcal{C}$
of real codimension 4 in $X$, and the holomorphic discs which intersect
$\mathcal{C}$ form a codimension 2 family. In
particular, if we choose the point $p\in L$ in the complement of a
codimension 2 subset of $L$ then none of the Maslov index 2 discs
whose boundary passes through $p$ hits $\mathcal{C}$. This allows us to
define $n_\beta(L)$. 

Similarly, in the presence of $J$-holomorphic spheres of negative Chern
number, there might exist stable maps in the class $\beta$ consisting of
a disc component of Maslov index $>2$ whose boundary passes through the
point $p$ together with multiply covered spheres of negative Chern number.
The moduli space of such maps typically has excess dimension. However,
suitable assumptions on spheres of negative Chern number ensure that
these stable maps cannot occur as limits of sequences of honest discs of
Maslov index 2 as long as $p$ stays away from a codimension 2 subset in $L$,
which allows us to ignore the issue.
\end{remark}

\begin{remark}\label{rmk:virtual}
In the above discussion we have avoided the use of virtual perturbation
techniques. However, at the cost of additional technical
complexity we can remove (2) and (3) from Assumption \ref{as:good}.
Indeed, even if holomorphic discs of Maslov index 2 fail to be regular, as long
as there are no holomorphic discs of Maslov index $\le 0$ we can still
define $n_\beta(L)$ as a {\it virtual} count. 
Namely, the minimality of the Maslov index prevents bubbling of discs, 
so that when $\mu(\beta)=2$ the virtual fundamental chain 
$[\mathcal{M}(L,\beta)]^{vir}$ is actually a cycle, and
$n_\beta(L)$ can be defined as the degree of the evaluation map.
Moreover, $n_\beta(L)$ is locally constant under Lagrangian isotopies
as long are discs of Maslov index $\le 0$ do not occur: indeed, the
Lagrangian isotopy induces a cobordism between the virtual fundamental
cycles of the moduli spaces. 
\end{remark}

\subsection{Maslov index zero discs and wall-crossing I}\label{ss:wall1}
In actual examples, condition (1) in Assumption \ref{as:good} almost
never holds (with the notable exception of the toric case). 
Generically, in dimension $n\ge 3$, the best we can hope for is:

\begin{assumption}
All simple (non multiply covered) nonconstant holomorphic discs of Maslov
index 0 in $(X,L)$ are regular, and the associated evaluation maps at
boundary marked points are transverse to each other and to the evaluation
maps at boundary marked points of holomorphic discs of Maslov index 2.
\end{assumption}

Then simple nonconstant holomorphic discs of Maslov index 0 occur
in $(n-3)$-dimensional families, and the set $\mathcal{Z}$ of points of
$L$ which lie on the boundary of a nonconstant Maslov index 0 disc has codimension 2 in $L$.
For a generic point $p\in L$, in each relative homotopy class of Maslov
index 2 there are finitely many holomorphic discs whose
boundary passes through $p$, and none of them hits $\mathcal{Z}$.
We can therefore define an integer $n_\beta(L,p)$ which counts these discs
with appropriate signs, and by summing over $\beta$ as in Definition
\ref{def:m0} we obtain a complex number $m_0(L,\nabla,p)$.

\begin{figure}[t]
\setlength{\unitlength}{5mm}
\begin{picture}(24,4)(1,1)
\psset{unit=\unitlength}
\newgray{ltgray}{0.85}
\pscustom[fillstyle=solid,fillcolor=ltgray]{
\psarc(2.99,3){2}{19}{341}
\psarcn(5.25,2.2){0.4}{160}{30}
\psarc(6.63,3){1.2}{220}{140}
\psarcn(5.25,3.8){0.4}{330}{200}}
\put(3,3.1){\makebox(0,0)[cb]{$\beta$}}
\put(3,2.7){\makebox(0,0)[ct]{\tiny($\mu=2$)}}
\put(1.268,2){\circle*{0.25}}
\put(1.1,1.85){\makebox(0,0)[rt]{$p$}}
\pscircle[fillstyle=solid,fillcolor=ltgray](12,3){2}
\pscircle[fillstyle=solid,fillcolor=ltgray](15.2,3){1.2}
\put(12,3.1){\makebox(0,0)[cb]{$\beta'$}}
\put(12,2.7){\makebox(0,0)[ct]{\tiny($\mu=2$)}}
\put(15.2,3.2){\makebox(0,0)[cb]{$\alpha$}}
\put(15.2,2.9){\makebox(0,0)[ct]{\tiny($\mu=0$)}}
\put(10.268,2){\circle*{0.25}}
\put(14,3){\circle*{0.25}}
\put(10.1,1.85){\makebox(0,0)[rt]{$p$}}
\put(10.8,1.2){\makebox(0,0)[rt]{\tiny $(p\in\mathcal{W})$}}
\pscircle[fillstyle=solid,fillcolor=ltgray](20.5,3){2}
\pscircle[linestyle=dashed, dash=1.5pt 1.5pt](24.7,3){1.2}
\put(20.5,3.1){\makebox(0,0)[cb]{$\beta'$}}
\put(20.5,2.7){\makebox(0,0)[ct]{\tiny($\mu=2$)}}
\put(24.7,2.9){\makebox(0,0)[cb]{$\alpha$}}
\put(18.768,2){\circle*{0.25}}
\put(22.5,3){\circle*{0.2}}
\put(23.5,3){\circle*{0.2}}
\put(18.6,1.85){\makebox(0,0)[rt]{$p$}}
\psline[linestyle=dashed, dash=1.5pt 1.5pt](22.5,3)(23.5,3)
\end{picture}
\caption{Wall-crossing for discs}
\label{fig:wallcross}
\end{figure}
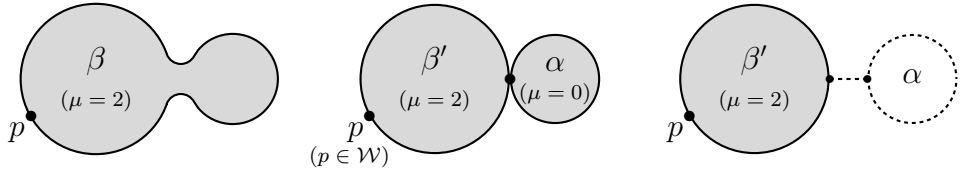

However, the points $p$ which lie on the boundary of a configuration
consisting of two holomorphic discs (of Maslov indices 2 and 0)
attached to each other at their boundary form a codimension 1 subset
$\mathcal{W}\subset L$. The typical behavior as $p$ approaches such a
``wall'' is that a Maslov index 2 disc representing a certain class $\beta$
breaks into a union of two discs representing classes $\beta'$ and
$\alpha$ with $\beta=\beta'+\alpha$, and then disappears altogether
(see Figure \ref{fig:wallcross}).
Thus the walls separate $L$ into various chambers,
each of which gives rise to a different value of $m_0(L,\nabla,p)$.

More conceptually,
denote by $\mathcal{M}_k(L,\beta)$ the moduli space of holomorphic discs
in $(X,L)$ with $k$ marked points on the boundary representing the class
$\beta$, and denote by $ev_i$ the evaluation map at the $i$-th marked point.
Then $n_\beta(L,p)$ is the degree at $p$ of the $n$-chain
$(ev_1)_*[\mathcal{M}_1(L,\beta)]$, whose boundary (an $n-1$-chain supported
on $\mathcal{W}$) is essentially
(ignoring all subtleties arising from multiple covers)
$$\sum_{\substack{\beta=\beta'+\alpha \\ \mu(\alpha)=0\\ 
0<\omega(\alpha)<\omega(\beta)}}
(ev_{1})_*[\mathcal{M}_2(L,\beta')\mathop{\times}\limits_{ev_2}
\mathcal{M}_1(L,\alpha)],$$
and $m_0(L,\nabla,p)$ is the degree at $p$ of the chain (with complex
coefficients) $$\m_0=\sum_\beta \exp(-\textstyle\int_\beta\omega)\,
\mathrm{hol}_\nabla(\partial\beta)\,(ev_1)_*[\mathcal{M}_1(L,\beta)].$$
In this language it is clear that these quantities depend on the position
of $p$ relatively to the boundary of the chain.

Various strategies can be employed to cancel the boundary and obtain
an evaluation cycle, thus leading to a well-defined count $n_\beta(L)$
independently of the point $p\in L$ \hbox{\cite{CL, FO3}}. For instance, in
the cluster approach \cite{CL}, given
a suitably chosen Morse function $f$ on $L$, one enlarges the moduli
space $\mathcal{M}_1(L,\beta)$ by considering configurations consisting
of several holomorphic discs connected to
each other by gradient flow trajectories of $f$, with one marked point
on the boundary of the component which lies at the root of the tree (which has Maslov index 2, while the
other components have Maslov index 0); see Figure \ref{fig:wallcross} (right)
for the simplest case.

However, even if one makes the disc count independent of the choice of 
$p\in L$ by completing the evaluation chain to a cycle, the final answer
still depends on the choice of auxiliary data.
For example, in the cluster construction, depending on the
direction of $\nabla f$ relative to the wall, two scenarios are possible:
either an honest disc in the class $\beta$ turns
into a configuration of two discs connected by a gradient flow line as
$p$ crosses $\mathcal{W}$; or
both configurations coexist on the same side of the wall (their contributions
to $n_\beta(L,p)$ cancel each other) and disappear as
$p$ moves across $\mathcal{W}$. Hence, in the absence of a canonical choice
there still isn't a uniquely defined superpotential.

\subsection{Maslov index zero discs and wall-crossing II: the surface case}
\label{ss:wall2}
The wall-crossing phenomenon is somewhat different in the surface case
($n=2$). In dimension 2 a generic Lagrangian submanifold does not bound any
holomorphic discs of Maslov index 0, so Assumption \ref{as:good} can be
expected to hold for most $L$, giving rise to a well-defined complex
number $m_0(L,\nabla)$. However, in a family of Lagrangians,
isolated holomorphic discs of Maslov index 0 occur in codimension 1,
leading to wall-crossing discontinuities. The general algebraic and
analytic framework which can be used to describe these phenomena is
discussed in \S 19.1 in \cite{FO3} (see also Section 5c in \cite{seidel}).
Here we discuss things in a more informal manner, in order to provide
some additional context for the calculations in Section \ref{s:nontoric}.

Consider a continuous family of (special) Lagrangian submanifolds $L_t$ $(t\in [-\epsilon,
\epsilon])$, such that $L_t$ satisfies Assumption \ref{as:good} for $t\neq 0$
and $L_0$ bounds a unique nontrivial simple holomorphic disc $u_\alpha$
representing a class $\alpha$ of Maslov index 0 (so
$\mathcal{M}(L_0,\alpha)=\{u_\alpha\}$). Given a holomorphic disc $u_0$
representing a class $\beta_0\in \pi_2(X,L_0)$ of Maslov index 2, we obtain
stable maps representing the class $\beta=\beta_0+m\alpha$ by attaching
copies of $u_\alpha$ (or branched covers of $u_\alpha$) to $u_0$ at
points where the boundary of $u_0$ intersects that of $u_\alpha$.
These configurations typically deform to honest holomorphic discs either
for $t>0$ or for $t<0$, but not both.

Using the isotopy to identify $\pi_2(X,L_t)$ with $\pi_2(X,L_0)$,
we can consider the moduli space of holomorphic discs with boundary
in one of the $L_t$, representing a given class $\beta$,
and with $k$ marked points on the boundary, $\tilde{\mathcal{M}}_k(\beta)=
\coprod_{t\in [-\epsilon,\epsilon]} \mathcal{M}_k(L_t,\beta)$, and
the evaluation maps $ev_i:\tilde{\mathcal{M}}_k(\beta)\to \coprod_t \{t\}\times L_t$.
In principle, given a class $\beta$ with $\mu(\beta)=2$, the boundary of
the corresponding evaluation chain is given by
\begin{equation}\label{eq:walld}
\partial\left((ev_1)_*[\tilde{\mathcal{M}}_k(\beta)]\right)=
\sum_{m\ge 1}  (ev_1)_*\biggl[\tilde{\mathcal{M}}_2(\beta-m\alpha)
\mathop{\times}\limits_{ev_2} \tilde{\mathcal{M}}_1(m\alpha)\biggr].
\end{equation}
However, interpreting the right-hand side of this equation is tricky,
because of the systematic failure of transversality, even if we ignore
the issue of multiply covered discs ($m\ge 2$). Partial relief can be
obtained by perturbing $J$ to a domain-dependent almost-complex structure.
Then, as one moves through the one-parameter family of
Lagrangians, bubbling of Maslov index 0 discs occurs at different values
of $t$ depending on the position at which it takes place along the boundary
of the Maslov index 2 component. So, as $t$ varies between
$-\epsilon$ and $+\epsilon$ one successively hits several boundary strata, 
corresponding to bubbling at various points of the boundary; the algebraic
number of such elementary wall-crossings is the intersection number
$[\partial\beta]\cdot[\partial\alpha]$.

As the perturbation of $J$ tends to zero, the various values of $t$ at which
wall-crossing occurs all tend to zero, and the representatives of the
classes $\beta-m\alpha$ which appear in the right-hand side of
(\ref{eq:walld}) might themselves undergo further bubbling as $t\to 0$.
Thus, we actually end up considering stable maps with boundary in $L_0$,
consisting of
several Maslov index 0 components (representing classes $m_i\alpha$)
attached simultaneously at $r$ different points of the boundary of a
main component of Maslov index 2. In a very informal sense, we can write
$$\vphantom{\bigg|}^{\mbox{``}}\,
\partial\left((ev_1)_*[\tilde{\mathcal{M}}_k(\beta)]\right)=\!\!\!\!\!\!\!\!
\sum_{\substack{m_1, \dots, m_r\ge 1\medskip\\
\beta=\beta'+(\sum m_i)\alpha}} \!\!\!\!\!\{0\}\times
(ev_1)_*\biggl[\mathcal{M}_{r+1}(L_0,\beta')
\mathop{\times}\limits_{ev_2,\dots,ev_{r+1}} \prod_{i=1}^r
\mathcal{M}_1(L_0,m_i\alpha)\biggr].
\vphantom{\bigg|}^{\mbox{''}}$$
However this formula is even more problematic than (\ref{eq:walld}), so
we will continue to use a domain-dependent almost-complex structure in
order to analyze wall-crossing.

On the other hand, one still has to deal with the failure of transversality
when the total contribution of the bubbles attached at a given point of
the boundary is a nontrivial multiple of $\alpha$.
Thus, in equation (\ref{eq:walld}) the moduli spaces associated to
multiple classes have to be unerstood in a virtual sense. Namely, for $m\ge
2$ we treat $\tilde{\mathcal{M}}(m\alpha)$ as a 0-chain
(supported at $t=0$) which corresponds
to the family count $\tilde{n}_{m\alpha}$ of discs in the class $m\alpha$
(e.g.\ after suitably perturbing the holomorphic curve equation). Typically,
for $m=1$ we have $\tilde{n}_\alpha=\pm 1$, while multiple cover contributions
are {\it a priori} harder to assess; in the example in \S \ref{s:nontoric}
they turn out to be zero, so since they should be determined by a purely local
calculation it seems reasonable to conjecture that they are always zero.
However at this point we do not care about the actual coefficients, all
that matters is that they only depend on the holomorphic disc of Maslov
index 0 ($u_\alpha$) and not on the class $\beta$.

Equation (\ref{eq:walld}) determines the manner in which the disc
counts $n_\beta(L_t)$ vary as $t$ crosses 0. It is easier to state
the formula in terms of a generating series which encodes the disc counts
in classes of the form $\beta_0+m\alpha$, namely 
$$F_t(q)=\sum_{m\in\Z} n_{\beta_0+m\alpha}(L_t)\,q^m.$$
Then each individual wall-crossing (at a given point on the boundary of the
Maslov index 2 disc) affects $F_t(q)$ by a same factor
$h_\alpha(q)=1+\tilde{n}_\alpha q+2\tilde{n}_{2\alpha} q^2+\dots$, so that in the end
$F_{-\epsilon}(q)$ and $F_{+\epsilon}(q)$ differ by a multiplicative factor of
$h_\alpha(q)^{[\partial \beta_0]\cdot[\partial\alpha]}$. 

Next, we observe that the contributions of the discs in the classes
$\beta_0+m\alpha$ to $m_0(L_t,\nabla_t)$ are given by plugging
$q=z_\alpha$ (as defined in (\ref{eq:coord})) into $F_t(q)$ and
multiplying by $z_{\beta_0}$. The values of this expression on either side
of $t=0$ differ from each other by a {\it change of variables},
replacing $z_{\beta_0}$ by
$z^*_{\beta_0}=z_{\beta_0}h_\alpha(z_\alpha)^{[\partial\beta_0]\cdot[\partial\alpha]}$.
These changes of variables can be performed consistently for all classes,
in the sense that the new variables still satisfy
$z^*_{\beta+\gamma}=z^*_\beta z^*_\gamma$. To summarize the discussion, we
have:

\begin{proposition}\label{prop:wallcross}
Upon crossing a wall in which $L$ bounds a unique simple Maslov index 0 
disc representing a relative class $\alpha$, the expression of
$m_0(L,\nabla)$ as a Laurent series in the variables
of Lemma \ref{l:holomcoords} is modified by a holomorphic change of
variables $$z_\beta\mapsto z_\beta \, h(z_\alpha)^{[\partial\beta]
\cdot[\partial\alpha]} \quad \forall \beta\in \pi_2(X,L),$$
where $h(z_\alpha)$ is a power series of the form $1+O(z_\alpha)$
(independent of $\beta$).
\end{proposition}

In view of Remark \ref{rmk:formal}, these properties also follow
formally from Fukaya-Oh-Ohta-Ono's construction of
$A_\infty$-homomorphisms associated to wall-crossing (Sections 19.1 and 30.9
of \cite{FO3}), as discussed by Seidel in Section 5c of \cite{seidel}.

An interesting consequence (especially in the light of the discussion in
\S \ref{s:m0c1}) is that, while the critical points of the superpotential
are affected by the wall-crossing, its critical values are not. 
Note however
that, since the change of variables can map a critical point to infinity
(see e.g.\ Section \ref{s:cp1cp1} for a family of special Lagrangian tori
on $\CP^1\times\CP^1$ in which this occurs), some critical values may still
be lost in the coordinate change.

Finally, we observe that the changes of variables which arise in
Proposition \ref{prop:wallcross} are formally very similar by the
quantum corrections to the complex structure of the mirror proposed by
Kontsevich-Soibelman and Gross-Siebert in the Calabi-Yau case
\cite{GS,KS2}. This suggests the following:

\begin{conj}\label{conj:quantumcorr}
The mirror to a K\"ahler surface $X$
(together with an anticanonical divisor $D$) should differ from
the complexified moduli space $M$ of special Lagrangian tori in $X\setminus D$
by ``quantum corrections'' which, away from the singular fibers, amount to
gluing the various regions of $M$ delimited by Maslov index 0 discs
according to the changes of variables introduced in Proposition
\ref{prop:wallcross}.
\end{conj}

One difficulty raised by this conjecture is that, whereas the quantum
corrections are compatible with the complex structure $J^\vee$, they do
not preserve the symplectic form $\omega^\vee$ introduced in Definition
\ref{def:omegavee}. We do not know how to address this issue, but
presumably this means that $\omega^\vee$ should also
be modified by quantum corrections.

\section{The toric case}\label{s:toric}

In this section, we consider the case where $X$ is a smooth toric variety,
and $D$ is the divisor consisting of all degenerate toric orbits. The
calculation of the superpotential (Proposition \ref{prop:toricW})
is very similar to that in \cite{ChoOh},
but we provide a self-contained description for completeness. 
We first recall very briefly some classical facts about toric varieties.

As a K\"ahler manifold, a toric variety $X$ is determined by its moment
polytope $\Delta\subset\R^n$, a convex polytope in which every
facet admits an integer normal vector, $n$ facets meet at
every vertex, and their primitive integer normal vectors form a basis of
$\Z^n$. The moment map $\phi:X\to\R^n$ identifies the orbit space of the
$T^n$-action on $X$ with $\Delta$. From the point of view of complex
geometry, the preimage of the interior of $\Delta$ is an open dense
subset $U$ of $X$, biholomorphic to $(\C^*)^n$, on which $T^n=(S^1)^n$ acts
in the standard manner. Moreover $X$ admits an open cover by affine
subsets biholomorphic to $\C^n$, which are the
preimages of the open stars of the vertices of $\Delta$ (i.e., the union
of all the strata whose closure contains the given vertex).

For each facet $F$ of $\Delta$, the preimage $\phi^{-1}(F)=D_F$ is a
hypersurface in $X$; the union of these hypersurfaces defines the toric
anticanonical divisor $D=\sum_F D_F$. The standard holomorphic volume 
form on $(\C^*)^n\simeq U=X\setminus D$, defined in coordinates by
$\Omega=d\log x_1\wedge\dots\wedge d\log x_n$, determines a section of
$K_X$ with poles along $D$.

\subsection{Toric orbits and the superpotential}
Our starting point is the observation that the moment map defines a
special Lagrangian torus fibration on $U=X\setminus D$:

\begin{lemma}
The $T^n$-orbits in $X\setminus D$ are special Lagrangian (with phase
$n\pi/2$).
\end{lemma}

\proof It is a classical fact that the $T^n$-orbits are Lagrangian; since
the $T^n$-action on $X\setminus D\simeq (\C^*)^n$ is the standard one, in
coordinates the orbits are products of circles $S^1(r_1)\times\dots\times
S^1(r_n)=\{(x_1,\dots,x_n),\ |x_i|=r_i\}$, on which the restriction of
$\Omega=d\log x_1\wedge\dots\wedge d\log x_n$ has constant phase $n\pi/2$.
\endproof

As above we consider the complexified moduli space $M$,
i.e.\ the set of pairs
$(L,\nabla)$ where $L$ is a $T^n$-orbit and $\nabla$ is a flat $U(1)$
connection on the trivial bundle over $L$. Recall that $L$ is a product
of circles $L=S^1(r_1)\times\dots\times S^1(r_n)\subset(\C^*)^n\simeq
X\setminus D$, and
denote the holonomy of $\nabla$ around the $j$-th factor $S^1(r_j)$ by
$\exp(i\theta_j)$. Then the symplectic form introduced in Definition
\ref{def:omegavee} becomes $\omega^\vee=(2\pi)^n
\sum d\log r_j\wedge d\theta_j$, i.e.\ up to a constant factor
it coincides with the standard K\"ahler form on $(\C^*)^n\simeq M$. 
However, as a complex manifold $M$ is not biholomorphic to $(\C^*)^n$:

\begin{proposition} \label{prop:Mtoric}
$M$ is biholomorphic to $\LLog^{-1}(\mathrm{int}\,\Delta)\subset (\C^*)^n$, where
$\LLog:(\C^*)^n\to\R^n$ is the map defined by 
$\LLog(z_1,\dots,z_n)=(-\frac{1}{2\pi}\log |z_1|,\dots,-\frac{1}{2\pi}\log |z_n|)$.
\end{proposition}

\proof
Given a $T^n$-orbit $L$ and a flat $U(1)$-connection $\nabla$, let
\begin{equation}\label{eq:torcoord}
z_j(L,\nabla)=\exp(-2\pi\phi_j(L))\,\hol_\nabla(\gamma_j),
\end{equation}
where $\phi_j$ is the $j$-th component of the moment map, i.e.\ the
Hamiltonian for the action of the $j$-th factor of $T^n$, and
$\gamma_j=[S^1(r_j)]\in H_1(L)$ is the homology class corresponding
to the $j$-th factor in $L=S^1(r_1)\times\dots\times S^1(r_n)$.

Let $A_j$ be a relative homology class in $H_2(X,L)$ such that
$\partial A_j=\gamma_j\in H_1(L)$ (it is clear that such a class can
be chosen consistently for all $T^n$-orbits), and consider the holomorphic
function $z_{A_j}$ defined by (\ref{eq:coord}): then $z_j$ and $z_{A_j}$
differ by a constant multiplicative factor. Indeed, comparing the two
definitions the holonomy factors coincide, and given an infinitesimal
special Lagrangian deformation $v\in C^\infty(NL)$,
$$\textstyle d\log|z_{A_j}|(v)=\int_{\gamma_j}-\iota_v\omega=\int_{S^1(r_j)}
\omega(X_j,v)\,dt=\int_{S^1(r_j)} -d\phi_j(v)\,dt=d\log |z_j|(v),$$
where $X_j$ is the vector field which generates the action of the
$j$-th factor of $T^n$ (i.e. the Hamiltonian vector field associated to
$\phi_j$).

Thus $z_1,\dots,z_n$ are holomorphic coordinates on $M$, and the
$(\C^*)^n$-valued map
$(L,\nabla)\mapsto (z_1(L,\nabla),\dots,z_n(L,\nabla))$ identifies $M$
with its image, which is exactly the preimage by $\LLog$ of the interior
of $\Delta$.
\endproof

Next we study holomorphic discs in $X$ with boundary on a given $T^n$-orbit
$L$. For each facet $F$ of $\Delta$, denote by $\nu(F)\in\Z^n$ the primitive
integer normal vector to $F$ pointing into $\Delta$, and let $\alpha(F)\in\R$
be the constant such that the equation of $F$ is $\langle \nu(F),\phi\rangle+
\alpha(F)=0$. Moreover, given $a=(a_1,\dots,a_n)\in \Z^n$ we denote by
$z^a$ the Laurent monomial $z_1^{a_1}\dots z_n^{a_n}$, where $z_1,\dots,z_n$
are the coordinates on $M$ defined by (\ref{eq:torcoord}). Then we have:

\begin{proposition}[Cho-Oh \cite{ChoOh}]\label{prop:toricW}
There are no holomorphic discs of Maslov index 0 in $(X,L)$, and the discs
of Maslov index 2 are all regular. Moreover, the superpotential is given
by the Laurent polynomial
\begin{equation}\label{eq:toricW}
W=m_0(L,\nabla)=\sum_{F\ \mathrm{facet}} e^{-2\pi\alpha(F)}\,z^{\nu(F)}.
\end{equation}
\end{proposition}

\proof
By Lemma \ref{l:maslov} and positivity of intersection,
Maslov index 0 discs do not intersect $D$, and
hence are contained in $X\setminus D\simeq (\C^*)^n$. However, since $L$
is a product of circles $S^1(r_i)=\{|x_i|=r_i\}$ inside $(\C^*)^n$, it
follows immediately from the maximum principle applied to $\log x_i$ that
$(\C^*)^n$ does not contain any non-constant holomorphic disc with boundary in $L$.

Next, we observe that a Maslov index 2 disc intersects $D$ at a single point, and in
particular it intersects only one of the components, say $D_F$ for some
facet $F$ of $\Delta$. 
We claim that for each facet $F$ there is a unique such disc
whose boundary passes through a given point $p=(x^0_1,\dots,
x^0_n)\in L\subset(\C^*)^n\simeq X\setminus D$; in terms of the components
$(\nu_1,\dots,\nu_n)$ of the normal vector $\nu(F)$, this disc can be
parametrized by the map 
\begin{equation}\label{eq:nudisc}
w\mapsto (w^{\nu_1} x^0_1,\dots,w^{\nu_n} x^0_n)
\end{equation}
(for $w\in D^2\setminus\{0\}$; the point $w=0$ corresponds to the
intersection with $D_F$).

To prove this claim, we work in an affine chart centered at a
vertex $v$ of $\Delta$ adjacent to the facet $F$. Denote by 
$\eta_1,\dots,\eta_n$ the basis of $\Z^n$ which consists of the
primitive integer vectors along the edges of
$\Delta$ passing through $v$, oriented away from $v$, and labelled in
such a way that $\eta_2,\dots,\eta_n$ are tangent to $F$.
Associate to each edge vector $\eta_i=(\eta_{i1},\dots,\eta_{in})\in\Z^n$ 
a Laurent monomial $\tilde{x}_i=x^{\eta_i}=x_1^{\eta_{i1}}\dots
x_n^{\eta_{in}}$. Then, after the change
of coordinates $(x_1,\dots,x_n)\mapsto (\tilde{x}_1,\dots,\tilde{x}_n)$,
the affine coordinate chart associated to the vertex $v$ can be thought of
as the standard compactification of $(\C^*)^n$ to $\C^n$. In this coordinate
chart, $L$ is again a product torus $S^1(\tilde{r}_1)\times\dots\times
S^1(\tilde{r}_n)$, where $\tilde{r}_i=r_1^{\eta_i1}\dots r_n^{\eta_{in}}$,
and $D_F$ is the coordinate hyperplane $\tilde{x}_1=0$. 

Since the complex
structure is the standard one, a holomorphic map $u:D^2\to\C^n$ with
boundary in $L$ is given by $n$ holomorphic functions
$w\mapsto (u_1(w),\dots,u_n(w))$ such that $|u_i|=\tilde{r}_i$ on the unit
circle. Since by assumption the disc hits only $D_F$, the functions
$u_2,\dots, u_n$ have no zeroes, so by the maximum principle
they are constant. Moreover the intersection number with $D_F$ is
assumed to be 1, so the image of the map $u_1$ is the disc of radius 
$\tilde{r}_1$, with multiplicity 1; so, up to reparametrization,
$u_1(w)=\tilde{r}_1\,w$. Thus, if we require the boundary of the disc
to pass through a given point $p=(\tilde{x}_1^0,\dots,\tilde{x}_n^0)$ of
$L$, then the only possible map (up to reparametrization) is
\begin{equation}\label{eq:trivdisc}
u:w\mapsto (w\,\tilde{x}_1^0,\tilde{x}_2^0,\dots,\tilde{x}_n^0),
\end{equation} 
which in the original coordinates is exactly (\ref{eq:nudisc}).

Moreover, it is easy to check (working component by component) that the map
(\ref{eq:trivdisc}) is regular. In particular, its contribution to the count of
holomorphic discs is $\pm 1$, and if we equip $L$ with the trivial spin
structure, then the sign depends only on the dimension $n$, and not on
the choice of the facet $F$ or of the $T^n$-orbit $L$. Careful inspection
of the sign conventions (see e.g.\ \cite{ChoOh,FO3,sebook}) shows that
the sign is $+1$.

The only remaining step in the proof of Proposition \ref{prop:toricW} is
to estimate the symplectic area of the holomorphic disc (\ref{eq:trivdisc}).
For this purpose, we first relabel the toric action so that it becomes
standard in the affine chart associated to the vertex $v$. Namely,
observe that the normal vectors $\nu(F_1)=\nu(F),\dots,\nu(F_n)$ to the
facets which pass through $v$ form a basis of $\Z^n$ dual to that given
by the edge vectors $\eta_1,\dots,\eta_n$. If we precompose the $T^n$-action
with the linear automorphism of $T^n$ induced by the transformation
$\sigma\in GL_n(\Z)$ which maps the $i$-th vector of the standard basis
to $\nu(F_i)$, then the relabelled action becomes standard in the
coordinates $(\tilde{x}_1,\dots,\tilde{x}_n)$. 

After this relabelling, the moment map becomes $\tilde{\phi}=
\sigma^T\circ \phi$, and in a neighborhood of the vertex 
$\tilde{v}=\sigma^T(v)$ the moment polytope $\tilde{\Delta}=\sigma^T(\Delta)$
is a translate of the standard octant. In particular, denoting by
$\tilde{\phi}_1$ the first component of $\tilde\phi$, the equation of the
facet $\tilde{F}=\sigma^T(F)$ is simply $\tilde{\phi}_1=-\alpha(F)$.
Since $u$ is equivariant with respect to the action of the first $S^1$ 
factor, integrating over the unit disc in polar coordinates $w=\rho
e^{i\theta}$ we have
$$\int_{D^2} u^*\omega=\iint_{D^2} \omega(\partial_\rho u,\partial_\theta u)
\,d\rho\,d\theta=\int_0^{2\pi}\int_0^1
d\tilde{\phi}_1(\partial_\rho u)\,d\rho\,d\theta=2\pi(\tilde\phi_1(L)-
\tilde\phi_1(u(0))).$$
Since $u(0)\in D_F$, we conclude that
$$\int_{D^2} u^*\omega=
2\pi(\tilde\phi_1(L)+\alpha(F))=2\pi\langle \nu(F),\phi(L)\rangle+
2\pi\alpha(F).$$
Incorporating the appropriate holonomy factor, we conclude that the
contribution of $u$ to the superpotential is precisely $e^{-2\pi\alpha(F)}z^{\nu(F)}$.
\endproof

\subsection{Comparison with the Hori-Vafa mirror and renormalization}
\label{ss:compare}
The formula (\ref{eq:toricW}) is identical to the well-known formula for
the superpotential of the mirror to a toric manifold (see Section 5.3 of
\cite{HV}). However, our mirror is ``smaller'' than the usual one, because
the variables $(z_1,\dots,z_n)$ are constrained to lie in a bounded subset
of $(\C^*)^n$. In particular, since the norm of each term in the sum
(\ref{eq:toricW}) is bounded by 1 (as the symplectic area of a holomorphic
disc is always positive), in our situation $W$ is always
bounded by the number of facets of the moment polytope $\Delta$.
While the ``usual'' mirror could be recovered by analytic continuation
from $M$ to all of $(\C^*)^n$ (or equivalently, by allowing the holonomy
of the flat connection to take values in $\C^*$ instead
of $U(1)$), there are various reasons for not proceeding in this manner,
one of them being that the symplectic form $\omega^\vee$ on $M$ blows up near
the boundary.

In fact, our description of $M$ resembles very closely one of the
intermediate steps in Hori and Vafa's construction (see Section
3.1 of \cite{HV}). The dictionary between the two constructions is the
following. Given a facet $F$ of $\Delta$, let
$y_F=2\pi\alpha(F)-\log(z^{\nu(F)})$, so that the real part of $y_F$ is
the symplectic area of one of the Maslov index 2 discs bounded by $L$
and its imaginary part is (minus) the integral of the connection 1-form
along its boundary. Then $M$ is precisely the subset of $(\C^*)^n$ in
which $\Re(y_F)>0$ for all facets, and the K\"ahler form $\omega^\vee$
introduced in Definition \ref{def:omegavee} blows up for $\Re(y_F)\to 0$.
This is exactly the same behavior as in equation (3.22) of \cite{HV}
(which deals with the case of a single variable $y$). 

Hori and Vafa introduce a renormalization procedure which enlarges
the mirror and flattens its K\"ahler metric. While the mathematical
justification for this procedure is somewhat unclear, it is interesting to
analyze it from the perspective of our construction. Hori and Vafa's
renormalization process replaces the inequality $\Re(y_F)>0$ by
$\Re(y_F)>-k$ for some constant $k$ (see equations (3.24) and (3.25)
in \cite{HV}), without changing the formula for the superpotential.
This amounts to enlarging the moment polytope by $\frac{1}{2\pi}k$
in all directions.

Assuming that $X$ is Fano (or more generally that $-K_X$ is nef),
another way to enlarge $M$ in the same manner is to equip $X$ with a
``renormalized'' K\"ahler form $\omega_k$ (compatible with the toric
action) chosen so that $[\omega_k]=[\omega]+k\,c_1(X)$.
Compared to Hori and Vafa's renormalization, this operation has the
effect of not only extending the domain of definition of the
superpotential, but also rescaling it by a factor of $e^{-k}$;
however, if we simultaneously rescale the K\"ahler form on $X$
and the superpotential, then we obtain a result consistent with Hori
and Vafa's. This suggests:

\begin{conj}\label{conj:renormalize}
The construction of the mirror to a Fano manifold $X$ should be carried
out not by using the fixed K\"ahler form $\omega$, but instead by considering
a family of K\"ahler forms in the classes $[\omega_k]=[\omega]+k\,c_1(X)$,
equipping the corresponding complexified moduli spaces of special Lagrangian
tori with the rescaled superpotentials $e^k\, W_{(\omega_k)}$, and taking the limit as
$k\to+\infty$.
\end{conj}

Of course, outside of the toric setting it is not clear what it means to
``take the limit as $k\to+\infty$''. A reasonable guess is that one
should perform {\it symplectic inflation} along $D$, i.e.\ modify
the K\"ahler form by (1,1)-forms supported in a small neighborhood
$\mathcal{V}$ of $D$, constructed e.g.\ as suitable smooth approximations
of the (1,1)-current dual to $D$. Special Lagrangians
which lie in the complement of $\mathcal{V}$ are not affected by this
process:  given $L\subset X\setminus\mathcal{V}$, the only
effect of the inflation procedure is that the symplectic areas of the
Maslov index 2 discs bounded by $L$ are increased by $k$; this is precisely
compensated by the rescaling of the superpotential by a multiplicative
factor of $e^k$. On the other hand, near
$D$ the change of K\"ahler form should ``enlarge'' the moduli space of special
Lagrangians.

In the non-Fano case (more specifically when $-K_X$ is not nef), it is
not clear how renormalization should be performed, or even whether
it should be performed at all. For example, consider the
Hirzebruch surface $\mathbb{F}_m=\mathbb{P}(\mathcal{O}_{\mathbb{P}^1}\oplus
\mathcal{O}_{\mathbb{P}^1}(m))$ for $m>2$, as studied in \S 5.2 of
\cite{AKO1}. The superpotential is given
by $$W=z_1+z_2+\frac{e^{-A}}{z_1z_2^m}+\frac{e^{-B}}{z_2},$$
where $A$ and $B$ are the symplectic areas of the zero section (of square
$+m$) and the fiber respectively, satisfying $A>mB$. An easy calculation
shows that $W$ has $m+2$ critical points in $(\C^*)^2$; 
the corresponding vanishing cycles generate the Fukaya
category of this Landau-Ginzburg model. As explained in \cite{AKO1},
this is incorrect from the point 
of view of homological mirror symmetry, since the derived category of coherent
sheaves of $\mathbb{F}_m$ is equivalent to a subcategory generated by
only four of these $m+2$ vanishing cycles.
An easy calculation shows that the $z_2$ coordinates of the critical
points are the roots of the equation
$$z_2^{m-2}(z_2^2-e^{-B})^2-m^2 e^{-A}=0.$$
Provided that $A>mB+O(\log m)$, one easily shows that only four of the
roots lie in the range $e^{-B}<|z_2|<1$ (and these satisfy
$|z_1|<1$ and $|e^{-A}/z_1z_2^{m}|<1$ as needed). This suggests that
one should only consider $M=\LLog^{-1}(\mathrm{int}\,\Delta)\subset(\C^*)^2$
rather than all of $(\C^*)^2$.
(Note however that the behavior is not quite the expected one when $A$
is too close to $mB$, for reasons that are not entirely clear).

Perhaps a better argument against renormalization (or analytic continuation)
in the non-Fano case can be found in Abouzaid's work \cite{Ab1,Ab2}.
Abouzaid's approach to homological mirror symmetry for toric varieties
is to consider admissible Lagrangians which occur as sections of the $\LLog$
map with boundary in a suitable tropical deformation of the fiber of the
Landau-Ginzburg model. More specifically, the deformed fiber lies near
$\LLog^{-1}(\Pi)$, where $\Pi$ is the tropical hypersurface in $\R^n$
associated to a rescaling of the Laurent polynomial $W$;
the interior of $\Delta$ is a connected component of $\R^n\setminus\Pi$,
and Abouzaid only considers admissible Lagrangian sections of
the $\LLog$ map over this connected component. Then the results in 
\cite{Ab1,Ab2} establish a correspondence between these
sections and holomorphic line bundles over $X$.

When $-K_X$ is not nef (for example, for the Hirzebruch
surface $\mathbb{F}_m$ with $m>2$), $\R^n\setminus\Pi$ has more
than one bounded connected component, and the other components also
give rise to some admissible Lagrangians; however Abouzaid's work shows
that those are not relevant to mirror symmetry for $X$, and that one
should instead focus exclusively on those Lagrangians which lie in the
bounded domain $M\subset(\C^*)^n$.

\section{A non-toric example}\label{s:nontoric}

The goal of this section is to work out a specific example in the
non-toric setting, in order to illustrate some general features which
are not present in the toric case, such as wall-crossing phenomena and quantum corrections.
Let $X=\CP^2$, equipped with the standard Fubini-Study K\"ahler form (or
a multiple of it), and consider the anticanonical divisor
$D=\{(x:y:z),\ (xy-\epsilon z^2)z=0\}$ (the union of the conic $xy=\epsilon
z^2$ and the line $z=0$), for some $\epsilon\neq 0$. We equip
$\CP^2\setminus D$ with the holomorphic (2,0)-form 
which in the affine coordinate chart $\{(x:y:1),\
(x,y)\in\C^2\}$ is given by $$\Omega=\frac{dx\wedge dy}{xy-\epsilon}.$$

\subsection{A family of special Lagrangian tori}\label{ss:nontoric1}
The starting point of our construction is the pencil of conics defined
by the rational map $f:(x:y:z)\mapsto
(xy:z^2)$. We will mostly work in affine coordinates, and think of $f$
as the map from $\C^2$ to $\C$ defined by $f(x,y)=xy$, suitably extended
to the compactification. The fiber of $f$ above any non-zero complex
number is a smooth conic, while the fiber over $0$ is the union of two
lines (the $x$ and $y$ coordinate axes), and the fiber over $\infty$ is a
double line.

The group $S^1$ acts on each fiber of $f$ by $(x,y)\mapsto (e^{i\theta}x,
e^{-i\theta}y)$. We will consider Lagrangian tori which are contained in
$f^{-1}(\gamma)$ for some simple closed curve $\gamma\subset\C$, and consist
of a single $S^1$-orbit inside each fiber a point of $\gamma$.
Recall that the symplectic fibration $f$ carries a natural horizontal
distribution, given at every point by the symplectic orthogonal to the
fiber. Parallel transport with respect to this horizontal distribution
yields symplectomorphisms between smooth fibers, and
$L\subset f^{-1}(\gamma)$ is Lagrangian if and only if it is
invariant by parallel transport along $\gamma$. 

Each fiber of $f$ is foliated by $S^1$-orbits, and contains
a distinguished orbit that we call the {\it equator}, namely the set of
points where $|x|=|y|$. We denote by $\delta(x,y)$ the
signed symplectic area of the region between the $S^1$-orbit through $(x,y)$
and the equator in the fiber $f^{-1}(xy)$, with the convention that
$\delta(x,y)$ is positive if $|x|>|y|$ and negative if $|x|<|y|$.
Since $S^1$ acts by symplectomorphisms, parallel transport is 
$S^1$-equivariant. Moreover, the symplectic involution $(x,y)\mapsto (y,x)$
also preserves the fibers of $f$, and so parallel transport commutes with it. 
This implies that parallel transport maps equators to equators, and maps
other $S^1$-orbits to $S^1$-orbits in a $\delta$-preserving manner.

\begin{definition}
Given a simple closed curve $\gamma\subset\C$ and a real number
$\lambda\in(-\Lambda,\Lambda)$ (where $\Lambda=\int_{\CP^1}\omega$ is the
area of a line), we define $$T_{\gamma,\lambda}=\{(x,y)\in f^{-1}(\gamma),
\ \delta(x,y)=\lambda\}.$$ By construction $T_{\gamma,\lambda}$ is an
embedded Lagrangian torus in $\CP^2$, except when $0\in\gamma$ and
$\lambda=0$ (in which case it has a nodal singularity at the origin).

Moreover, when $0\not\in\gamma$, we say that $T_{\gamma,\lambda}$ is of {\em 
Clifford type} if $\gamma$ encloses the origin, and of {\em Chekanov
type} otherwise.
\end{definition}

This terminology is motivated by the observation that the product tori
$S^1(r_1)\times S^1(r_2)\subset\C^2$ (among which the Clifford
tori) are of the form $T_{\gamma,\lambda}$ where $\gamma$ is the circle
of radius $r_1r_2$ centered at the origin, whereas one way to define the
so-called {\it Chekanov torus} \cite{chekanov,EP} is as $T_{\gamma,0}$
for $\gamma$ a loop that does not enclose the origin (see \cite{EP}).

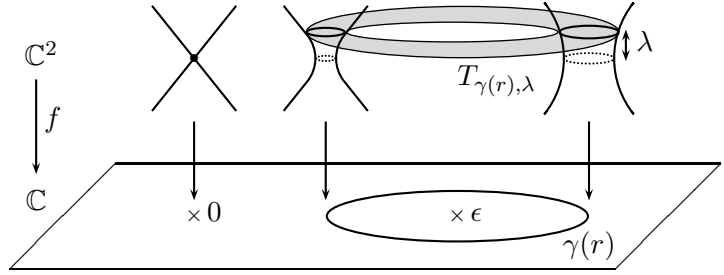
\begin{figure}[t]
\setlength{\unitlength}{7mm}
\begin{picture}(13.5,5)(0,0)
\psset{unit=\unitlength}
\newgray{ltgray}{0.85}
\psellipse[linewidth=0.5pt,fillstyle=solid,fillcolor=ltgray](8.6,4.5)(2.98,0.5)
\psellipse[linewidth=0.5pt,fillstyle=solid,fillcolor=white](8.4,4.5)(2.02,0.2)
\put(0,0){\line(1,0){11.5}}
\put(0,0){\line(1,1){2}}
\put(2,2){\line(1,0){11.5}}
\put(11.5,0){\line(1,1){2}}
\put(3.5,1){\makebox(0,0)[cc]{\tiny$\times$}}
\put(3.75,0.9){\small $0$}
\put(8.5,1){\makebox(0,0)[cc]{\tiny$\times$}}
\put(8.75,0.9){\small $\epsilon$}
\put(10.5,0.3){\small $\gamma(r)$}
\psellipse(8.5,1)(2.5,0.5)
\psline[linearc=0](2.7,3)(3.5,4)(2.7,5)
\psline[linearc=0](4.3,3)(3.5,4)(4.3,5)
\psline[linearc=0.7](5.2,3)(6,4)(5.2,5)
\psline[linearc=0.7](6.8,3)(6,4)(6.8,5)
\psline[linearc=1.7](10.2,3)(11,4)(10.2,5)
\psline[linearc=1.7](11.8,3)(11,4)(11.8,5)
\put(3.5,4){\circle*{0.15}}
\psellipse[linestyle=dotted,dotsep=0.6pt](6,4)(0.2,0.07)
\psellipse[linestyle=dotted,dotsep=0.6pt](11,4)(0.5,0.12)
\psellipse(11,4.5)(0.58,0.14)
\psellipse(6,4.5)(0.38,0.1)
\psline{<->}(11.7,3.95)(11.7,4.55)
\put(11.9,4.15){\small $\lambda$}
\psline{->}(3.5,2.8)(3.5,1.3)
\psline{->}(6,2.8)(6,1.3)
\psline{->}(11,2.8)(11,1.3)
\put(0.5,1.3){\makebox(0,0)[cc]{\small $\C$}}
\put(0.5,4.1){\makebox(0,0)[cc]{\small $\C^2$\!\!}}
\put(0.65,2.7){\small $f$}
\psline{->}(0.5,3.6)(0.5,1.8)
\put(8.5,3.5){\small $T_{\gamma(r),\lambda}$}
\end{picture}
\caption{The special Lagrangian torus $T_{\gamma(r),\lambda}$} \label{fig:tcp2}
\end{figure}

Recall that the anticanonical divisor $D$ is the union of the fiber
$f^{-1}(\epsilon)$ and the line at infinity. The following proposition
motivates our interest in the tori $T_{\gamma,\lambda}$ in the specific
case where $\gamma=\gamma(r)$ is a circle
of radius $r$ centered at $\epsilon$. 

\begin{proposition}\label{prop:cp2slag}
The tori $T_{\gamma(r),\lambda}=\{(x,y),\ |xy-\epsilon|=r,\
\delta(x,y)=\lambda\}$ are special Lagrangian with respect to $\Omega=
(xy-\epsilon)^{-1}\,dx\wedge dy$.
\end{proposition}

\proof
Let $H(x,y)=|xy-\epsilon|^2$, and let $X_H$ be the corresponding
Hamiltonian vector field, i.e.\ the vector field such that
$\iota_{X_H}\omega=dH$. We claim that $X_H$ is everywhere tangent to
$T_{\gamma(r),\lambda}$. In fact, $H$ is constant over each fiber of $f$,
so $X_H$ is symplectically orthogonal to the fibers, i.e.\ it lies in the
horizontal distribution. Moreover, $X_H$ is tangent to the level sets of
$H$; so, up to a scalar factor, $X_H$ is in fact the horizontal lift of
the tangent vector to $\gamma(r)$, and thus it is tangent to $T_{\gamma(r),
\lambda}$. 

The tangent space to $T_{\gamma(r),\lambda}$ is therefore
spanned by $X_H$ and by the vector field generating the $S^1$-action,
$\xi=(ix,-iy)$. However, we observe that $$\iota_\xi\Omega=
\frac{ix\,dy+iy\,dx}{xy-\epsilon}=i\,d\log(xy-\epsilon).$$
It follows that $\Im\Omega(\xi,X_H)=d\log|xy-\epsilon|\,(X_H)$, which
vanishes since $X_H$ is tangent to the level sets of $H$. Hence
$T_{\gamma(r),\lambda}$ is special Lagrangian.
\endproof

Thus $\CP^2\setminus D$ admits a fibration by special Lagrangian tori
$T_{\gamma(r),\lambda}$, with a single nodal fiber $T_{\gamma(|\epsilon|),
0}$. For $r<|\epsilon|$ the tori $T_{\gamma(r),\lambda}$ are of Chekanov type,
while for $r>|\epsilon|$ they are of Clifford type. We shall now see that
wall-crossing occurs for $r=|\epsilon|$, thus separating the moduli space
into two chambers $r<|\epsilon|$ and $r>|\epsilon|$. We state the next two
lemmas in a slightly more general context.

\begin{lemma} \label{l:cp2maslov}
If $\gamma\subset\C$ is a simple closed loop and $w\in\C$ lies
in the interior of $\gamma$, then for any class $\beta\in\pi_2(\CP^2,
T_{\gamma,\lambda})$, the Maslov index is $\mu(\beta)=2(\beta\cdot [f^{-1}(w)]+
\beta\cdot [\CP^1_\infty])$, where $\CP^1_\infty$ is the line at infinity
in $\CP^2$.
\end{lemma}

\proof
If $\gamma$ is a circle centered at $w$, then Proposition
\ref{prop:cp2slag} implies that $T_{\gamma,\lambda}$ is special Lagrangian
for $\Omega=(xy-w)^{-1}dx\wedge dy$, and the result is then a direct
consequence of Lemma \ref{l:maslov}. The general case follows by
continuously deforming $\gamma$ to such a circle, without crossing $w$
nor the origin, and keeping track of relative homotopy classes through this
Lagrangian deformation, which affects neither the Maslov index nor the
intersection numbers with $f^{-1}(w)$ and $\CP^1_\infty$.
\endproof

Using positivity of intersection, this lemma precludes the existence of
holomorphic discs with negative Maslov index. Moreover:

\begin{lemma} \label{l:cp2wall}
The Lagrangian torus $T_{\gamma,\lambda}$ bounds a nontrivial holomorphic
disc of Maslov index 0 if and only if\/ $0\in\gamma$.
\end{lemma}

\proof
Assume there is a non-trivial holomorphic map $u:(D^2,\partial D^2)
\to (\CP^2,T_{\gamma,\lambda})$ representing a class of Maslov index 0,
and choose a point $w\in\C$ inside the region delimited by $\gamma$. By positivity of 
intersection and Lemma \ref{l:cp2maslov}, the image of $u$
must be disjoint from $f^{-1}(w)$ and from the line at infinity.
The projection $f\circ u$ is therefore a well-defined holomorphic map
from $(D^2,\partial D^2)$ to $(\C,\gamma)$, whose image avoids $w$.

It follows that $f\circ u$ is constant, i.e.\ the image of $u$ is contained
in the affine part of a fiber of $f$, say $f^{-1}(c)$ for some $c\in\gamma$.
However, for $c\neq 0$ the affine conic $xy=c$ is topologically a cylinder
$S^1\times\R$, intersected by $T_{\gamma,\lambda}$ in an essential circle,
which does not bound any nontrivial holomorphic disc. Therefore $c=0$,
and $0\in\gamma$.

Conversely, if $0\in \gamma$, we observe that $f^{-1}(0)$ is the union of
two complex lines (the $x$ and $y$ coordinate axes), and its intersection
with $T_{\gamma,\lambda}$ is a circle in one of them (depending on the
sign of $\lambda$). Excluding the degenerate case $\lambda=0$, it follows
that $T_{\gamma,\lambda}$ bounds a holomorphic disc of area $|\lambda|$,
contained in one of the coordinate axes; by Lemma \ref{l:cp2maslov}
its Maslov index is 0.
\endproof

\subsection{The superpotential}
We now consider the complexified moduli space $M$ associated to the
family of special Lagrangian tori constructed in Proposition
\ref{prop:cp2slag}. The goal of this section is to compute the
superpotential; by Lemma \ref{l:cp2wall}, the cases $r<|\epsilon|$ and
$r>|\epsilon|$ should be treated separately. 

We start with the Clifford case ($r>|\epsilon|$). By deforming
continuously $\gamma(r)$ into a circle centered at the origin without
crossing the origin, we obtain a Lagrangian isotopy from 
$T_{\gamma(r),\lambda}$ to a product torus $S^1(r_1)\times S^1(r_2)
\subset\C^2$, with the property that the minimal Maslov index of a
holomorphic disc remains at least 2 throughout the deformation.
Therefore, by Remark \ref{rmk:virtual}, for each class $\beta$
of Maslov index 2, the disc count $n_\beta(L)$ remains constant throughout
the deformation. The product torus corresponds to the toric case
considered in Section \ref{s:toric}, so we can use Proposition \ref{prop:toricW}.

Denote by $z_1$ and $z_2$ respectively the holomorphic coordinates
associated to the relative homotopy classes $\beta_1$ and $\beta_2$
of discs parallel to the $x$ and $y$ coordinate
axes in $(\C^2,S^1(r_1)\times S^1(r_2))$ via the formula
(\ref{eq:coord}). Then Proposition \ref{prop:toricW} implies:

\begin{proposition}\label{prop:cliffordW}
For $r>|\epsilon|$, the superpotential is given by
\begin{equation}\label{eq:cliffordW}
W=z_1+z_2+\dfrac{e^{-\Lambda}}{z_1z_2}.\end{equation}
\end{proposition}

The first two terms in this expression 
correspond to sections of $f$ over the disc $\Delta$ of radius $r$
centered at $\epsilon$ (the first one intersecting $f^{-1}(0)$ at a point
of the $y$-axis, while the second one hits the $x$-axis), whereas the
last term corresponds to a disc whose image under $f$ is a double cover
of $\CP^1\setminus \Delta$ branched at infinity.
\bigskip

Next we consider the case $r<|\epsilon|$, where $\gamma=\gamma(r)$ does not
enclose the origin. We start with the special case $\lambda=0$, which is
the one considered by Chekanov and Eliashberg-Polterovich \cite{chekanov,EP}. 

The fibration $f$ is trivial over the disc $\Delta$ bounded by
$\gamma$, and over $\Delta$ it admits an obvious holomorphic section with boundary in
$T_{\gamma,0}$, given by the portion of the line $y=x$ for which $x\in
\sqrt{\Delta}$ (one of the two preimages of $\Delta$ under $z\mapsto z^2$).
More generally, by considering the portion of the line
$y=e^{2i\theta}x$ where $x\in e^{-i\theta}\sqrt{\Delta}$, and letting
$e^{i\theta}$ vary in $S^1$, we obtain a family of holomorphic discs of
Maslov index 2 with boundary in $T_{\gamma,0}$. One easily checks that
these discs are regular, and that they boundaries sweep out $T_{\gamma,0}$
precisely once; we denote their class by $\beta$.

Other families of Maslov index 2 discs are harder to come by; the
construction of one such family is outlined in an unfinished manuscript
of Blechman and Polterovich, but the complete classification has only
been carried out recently by Chekanov and Schlenk \cite{schlenk}.
In order to state Chekanov and Schlenk's results, we need one more piece
of notation. Given a line segment which joins the origin
to a point $c=\rho e^{i\theta}\in\gamma$, consider the Lefschetz thimble associated to the
critical point of $f$ at the origin, i.e.\ the Lagrangian disc with
boundary in $T_{\gamma,0}$ formed by the collection of equators in the 
fibers of $f$ above the segment $[0,c]$; this is just
a disc of radius $\sqrt{\rho}$ in the line $y=e^{i\theta}\bar{x}$. We denote by $\alpha\in\pi_2(\CP^2,
T_{\gamma,0})$ the class of this disc; one easily checks that $\alpha$,
$\beta$, and $H=[\CP^1]$ form a basis of $\pi_2(\CP^2,T_{\gamma,0})$.

\begin{lemma}[Chekanov-Schlenk \cite{schlenk}]
The only classes in $\pi_2(\CP^2,T_{\gamma,0})$ which may contain
holomorphic discs of Maslov index 2 are $\beta$\/ and\/
$H-2\beta+k\alpha$ for $k\in\{-1,0,1\}$.
\end{lemma}

\proof
We compute the intersection numbers of $\alpha$, $\beta$ and
$H$ with the $x$-axis, the $y$-axis, and the fiber $f^{-1}(\epsilon)$,
as well as their Maslov indices (using Lemma \ref{l:cp2maslov}):
\begin{center}
\begin{tabular}{|c||c|c|c||c|}
\hline
class&$x$-axis&$y$-axis&$f^{-1}(\epsilon)\vphantom{\Big|}$&\ $\mu$\ \ \\
\hline
$\alpha\vphantom{\Big|}$&\!\!$-1$&1&0&0\\
$\beta$&\,\,0&0&1&2\\
$H\vphantom{\Big|}$&\,\,1&1&2&6\\
\hline
\end{tabular}
\end{center}
A class of Maslov index 2 is of the form $\beta+m(H-3\beta)+k\alpha$
for $m,k\in\Z$; the constraints on $m$ and $k$ come from positivity
of intersections. Considering the intersection number with
$f^{-1}(\epsilon)$, we must have $m\le 1$; and considering the intersection
numbers with the $x$-axis and the $y$-axis, we must have $m\ge |k|$. It
follows that the only possibilities are $m=k=0$ and $m=1$, $|k|\le 1$.
\endproof

\begin{proposition}[Chekanov-Schlenk \cite{schlenk}]
The torus $T_{\gamma,0}$ bounds a unique $S^1$-family of holomorphic discs 
in each of the classes $\beta$ and $H-2\beta+k\alpha$,
$k\in\{-1,0,1\}$. These discs are regular, and the corresponding evaluation
maps have degree 2 for $H-2\beta$ and 1 for the other classes.
\end{proposition}

\proof[Sketch of proof]
We only outline the construction of the holomorphic discs in the classes
$H-2\beta+k\alpha$, following Blechman-Polterovich and
Chekanov-Schlenk. The reader is referred to \cite{schlenk} for details
and for proofs of uniqueness and regularity.

Let $\varphi$ be a biholomorphism from the unit disc $D^2$ to the complement
$\CP^1\setminus\Delta$ of the disc bounded by $\gamma$, parametrized so that
$\varphi(0)=\infty$ and $\varphi(a^2)=0$ for some $a\in(0,1)$, and
consider the double branched cover $\psi(z)=\varphi(z^2)$, which has a
pole of order 2 at the origin and simple roots at $\pm a$. We will construct
holomorphic maps $u:(D^2,\partial D^2)\to (\CP^2,T_{\gamma,0})$ such that
$f\circ u=\psi$. Let $$\tau_a(z)=\frac{z-a}{1-az},\quad
\tau_{-a}(z)=\frac{z+a}{1+az},\quad \mathrm{and}\quad
g(z)=\frac{z^2\,\psi(z)}{\tau_a(z)\,\tau_{-a}(z)}.$$ 
Since $\tau_{\pm a}$ are biholomorphisms of the unit disc mapping
$\pm a$ to $0$, the map $g$ is a nonvanishing holomorphic
function over the unit disc, and hence we can choose a square root
$\sqrt{g}$. Then for any $e^{i\theta}\in S^1$ we can consider the holomorphic maps
\begin{eqnarray}
z\mapsto \bigl(e^{i\theta}\,\tau_a(z)\,\tau_{-a}(z)\,\sqrt{g(z)}:
e^{-i\theta}\,\sqrt{g(z)}:z\bigr),\label{eq:disc1}\\
z\mapsto \bigl(e^{i\theta}\,\tau_a(z)\,\sqrt{g(z)}:
e^{-i\theta}\,\tau_{-a}(z)\,\sqrt{g(z)}:z\bigr),\label{eq:disc2}\\
z\mapsto \bigl(e^{i\theta}\,\sqrt{g(z)}:
e^{-i\theta}\,\tau_a(z)\,\tau_{-a}(z)\,\sqrt{g(z)}:z\bigr).\label{eq:disc3}
\end{eqnarray}
Letting $u$ be any of these maps, it is easy to check that
$f\circ u=\psi$, and that the first two components of $u$ have equal
norms when $|z|=1$ (using the fact that $|\tau_a(z)|=|\tau_{-a}(z)|=1$
for $|z|=1$). So in all cases $\partial D^2$ is mapped to $T_{\gamma,0}$.
One easily checks (e.g.\ using intersection numbers with the coordinate
axes) that the classes represented by these maps are 
$H-2\beta+k\alpha$ with $k=1,0,-1$ respectively for
(\ref{eq:disc1})--(\ref{eq:disc3}).

Chekanov and Schlenk show that these maps are regular, and that this list
is exhaustive \cite{schlenk}.
(In fact, since they enumerate discs whose boundary passes
through a given point of $T_{\gamma,0}$, they also introduce a fourth map which
differs from (\ref{eq:disc2}) by swapping $\tau_a$ and $\tau_{-a}$; however
this is equivalent to reparametrizing by $z\mapsto -z$).

Finally, the degrees of the evaluation maps are easily determined by
counting the number of values of $e^{i\theta}$ for which the
boundary of $u$ passes through a given point of $T_{\gamma,0}$; however it
is important to note here that, for the maps (\ref{eq:disc1}) and
(\ref{eq:disc3}), replacing $\theta$ by $\theta+\pi$ yields the same disc
up to a reparametrization ($z\mapsto -z$).
\endproof

By Lemma \ref{l:cp2wall} and Remark \ref{rmk:virtual}, the disc counts
remain the same in the general case (no longer assuming $\lambda=0$), 
since deforming $\lambda$ to $0$ yields a Lagrangian isotopy from
$T_{\gamma,\lambda}$ to $T_{\gamma,0}$ in the complement of $f^{-1}(0)$.
Therefore, denoting by $u$ and $w$ the holomorphic coordinates on $M$
associated to the classes $\beta$ and $\alpha$ respectively, we have:

\begin{proposition}\label{prop:chekanovW}
For $r<|\epsilon|$, the superpotential is given by 
\begin{equation}\label{eq:chekanovW}
W=u+\frac{e^{-\Lambda}}{u^2w}+2\,\frac{e^{-\Lambda}}{u^2}+
\frac{e^{-\Lambda}w}{u^2}=u+\frac{e^{-\Lambda}(1+w)^2}{u^2w}.
\end{equation}
\end{proposition}

\subsection{Wall-crossing, quantum corrections and monodromy}\label{ss:nontoric3}
In this section, we compare the two formulas obtained for the superpotential
in the Clifford and Chekanov cases (Propositions \ref{prop:cliffordW} and
\ref{prop:chekanovW}), in terms of wall-crossing at $r=|\epsilon|$.
We start with a simple observation:

\begin{lemma}\label{l:cp2compare}
The expressions (\ref{eq:cliffordW}) and (\ref{eq:chekanovW}) are related
by the change of variables $u=z_1+z_2$, $w=z_1/z_2$.
\end{lemma}

To see how this fits with the general discussion of wall-crossing
in \S \ref{ss:wall2}
and Proposition \ref{prop:wallcross}, we consider separately the two cases
$\lambda>0$ and $\lambda<0$. We use the same notations as in the previous
section concerning relative homotopy classes ($\beta_1,\beta_2$ on the
Clifford side, $\beta,\alpha$ on the Chekanov side) and the corresponding
holomorphic coordinates on $M$ ($z_1,z_2$ and $u,w$).

First we consider the case where $\lambda>0$, i.e.\ $T_{\gamma(r),\lambda}$
lies in the region where $|x|>|y|$. When $r=|\epsilon|$, $T_{\gamma(r),
\lambda}$ intersects the $x$-axis in a circle, which bounds a disc $u_0$ of
Maslov index 0. In terms of the basis used on the Clifford side, the class
of this disc is  $\beta_1-\beta_2$; on the Chekanov side it is $\alpha$. 

As $r$ decreases through $|\epsilon|$, two of the families of
Maslov index 2 discs discussed in the previous section survive the
wall-crossing: namely the family of holomorphic discs in the class
$\beta_2$ on the Clifford side becomes the family of discs in the class
$\beta$ on the Chekanov side, and the discs in the class $H-\beta_1-
\beta_2$ on the Clifford side become the discs in the class
$H-2\beta-\alpha$ on the Chekanov side. This correspondence between
relative homotopy classes determines the change of variables between
the coordinate systems $(z_1,z_2)$ and $(u,w)$ of the two charts on $M$
along the $\lambda>0$ part of the wall:
\begin{equation} \label{eq:classicalpos}
\left\{
\begin{array}{rclrcl}
\alpha&\leftrightarrow&\beta_1-\beta_2&w&\leftrightarrow&z_1/z_2\\
\beta&\leftrightarrow&\beta_2&u&\leftrightarrow&z_2\\
\,H-2\beta-\alpha&\leftrightarrow&H-\beta_1-\beta_2\qquad\quad&
e^{-\Lambda}/u^2w&\leftrightarrow&e^{-\Lambda}/z_1z_2
\end{array}
\right.
\end{equation}

However, with this ``classical'' interpretation of the geometry of $M$
the formulas (\ref{eq:cliffordW}) and (\ref{eq:chekanovW}) do not match
up, and the superpotential presents a wall-crossing discontinuity,
corresponding to the contributions of the various families of discs
that exist only on one side of the wall. As $r$ decreases through
$|\epsilon|$, holomorphic discs in the class $\beta_1$ break into the
union of a disc in the class $\beta=\beta_2$ and the exceptional disc
$u_0$, and then disappear entirely. Conversely, new discs in the classes
$H-2\beta$ and $H-2\beta+\alpha$ are generated by attaching $u_0$ to
a disc in the class $H-\beta_1-\beta_2=H-2\beta-\alpha$ at one or both
of the points where their boundaries intersect. Thus the correspondence
between the two coordinate charts across the wall should be corrected to:
\begin{equation} \label{eq:quantumpos}
\left\{
\begin{array}{rclrcl}\vphantom{\Big|}
\beta&\leftrightarrow&\{\beta_1,\beta_2\}&u&\leftrightarrow&z_1+z_2\\
\,H-2\beta+\{-1,0,1\}\alpha&\leftrightarrow&H-\beta_1-\beta_2\qquad&
\dfrac{e^{-\Lambda}(1+w)^2}{u^2w}&\leftrightarrow&\dfrac{e^{-\Lambda}}{z_1z_2}
\end{array}
\right.
\end{equation}
This corresponds to the change of variables $u=z_1+z_2$, $w=z_1/z_2$ as
suggested by Lemma \ref{l:cp2compare}; the formula for $w$ is the same
as in (\ref{eq:classicalpos}), but the formula for $u$ is affected by a
multiplicative factor $1+w$, from $u=z_2$ to $u=z_1+z_2=(1+w)z_2$. This
is precisely the expected behavior in view of Proposition \ref{prop:wallcross}.

\begin{remark}
Given $c\in\gamma$, the class
$\alpha=\beta_1-\beta_2\in\pi_2(\CP^2,T_{\gamma,\lambda})$ can be represented
by taking the portion of $f^{-1}(c)$ lying between
$T_{\gamma,\lambda}$ and the equator, which has symplectic area $\lambda$,
together with a Lagrangian thimble.
Therefore $|w|=\exp(-\lambda)$. In particular, for $\lambda\gg 0$ the
correction factor $1+w$ is $1+o(1)$.
\end{remark}

The case $\lambda<0$ can be analyzed in the same manner. For $r=|\epsilon|$ the
Lagrangian torus $T_{\gamma(r),\lambda}$ now intersects the $y$-axis in a
circle; this yields a disc of Maslov index 0 representing the class
$\beta_2-\beta_1=-\alpha$. The two families
of holomorphic discs that survive the wall-crossing are those in the classes
$\beta_1$ and $H-\beta_1-\beta_2$ on the Clifford side, which become
$\beta$ and $H-2\beta+\alpha$ on the Chekanov side. Thus, the coordinate
change along the $\lambda<0$ part of the wall is
\begin{equation} \label{eq:classicalneg}
\left\{
\begin{array}{rclrcl}
-\alpha&\leftrightarrow&\beta_2-\beta_1&w^{-1}&\leftrightarrow&z_2/z_1\\
\beta&\leftrightarrow&\beta_1&u&\leftrightarrow&z_1\\
\,H-2\beta+\alpha&\leftrightarrow&H-\beta_1-\beta_2\qquad\quad&
e^{-\Lambda}w/u^2&\leftrightarrow&e^{-\Lambda}/z_1z_2
\end{array}
\right.
\end{equation}
However, taking wall-crossing phenomena into account, the
correspondence should be modified in the same manner as above,
from (\ref{eq:classicalneg}) to (\ref{eq:quantumpos}), which again leads to
the change of variables
$u=z_1+z_2$, $w=z_1/z_2$; this time, the formula for $u$ is corrected by a
multiplicative factor $1+w^{-1}$, from $u=z_1$ to $u=z_1+z_2=(1+w^{-1})z_1$.

\begin{remark}
The discrepancy between the gluing formulas (\ref{eq:classicalpos}) and
(\ref{eq:classicalneg}) is due to the monodromy of the family of special
Lagrangian tori $T_{\gamma(r),\lambda}$ around the nodal fiber
$T_{\gamma(|\epsilon|),0}$. The vanishing cycle of the nodal degeneration
is the loop $\partial\alpha$, and in terms of the basis $(\partial \alpha,
\partial \beta)$ of $H_1(T_{\gamma(r),\lambda},\Z)$, the monodromy is the
Dehn twist
$$\begin{pmatrix}1&1\\0&1\end{pmatrix}.$$ This induces monodromy
in the affine structures on the moduli space
$\mathcal{B}=\{(r,\lambda)\}$ and its complexification $M$. Namely,
$M$ carries an integral (complex) affine structure given by the coordinates
$(\log z_1,\log z_2)$ on the Clifford chamber $|r|>\epsilon$ and the
coordinates $(\log w,\log u)$ on the Chekanov chamber $|r|<\epsilon$.
Combining (\ref{eq:classicalpos}) and (\ref{eq:classicalneg}), moving
around $(r,\lambda)=(|\epsilon|,0)$ induces the transformation $(w,u)\mapsto
(w,uw)$, i.e.\ $$(\log w,\log u)\mapsto
(\log w,\log u+\log w).$$ Therefore, in terms of the basis $(\partial_{\log u},
\partial_{\log w})$ of $TM$, the monodromy is given by the transpose matrix
$$\begin{pmatrix}1&0\\1&1\end{pmatrix}.$$
Taking quantum corrections into account, the discrepancy in 
the coordinate transformation formulas disappears
(the gluing map becomes (\ref{eq:quantumpos}) for
both signs of $\lambda$), but the monodromy remains the same.
Indeed, the extra factors brought in by the quantum
corrections, $1+w$ for $\lambda>0$ and $1+w^{-1}$ for $\lambda<0$, are
both of the form $1+o(1)$ for $|\lambda|\gg 0$.
\end{remark}

\subsection{Another example: $\CP^1\times\CP^1$}\label{s:cp1cp1}
We now briefly discuss a related example: consider $X=\CP^1\times\CP^1$,
equipped with a product K\"ahler form such that the two factors have equal areas,
let $D$ be the union of the two lines at infinity and the conic
$\{xy=\epsilon\}\subset\C^2$, and consider
the 2-form $\Omega=(xy-\epsilon)^{-1}\,dx\wedge dy$ on $X\setminus D$.

The main geometric features remain the same as before, the main difference
being that the fiber at infinity of $f:(x,y)\mapsto xy$ is now a union of
two lines $L^1_\infty=\CP^1\times\{\infty\}$ and $L^2_\infty=\{\infty\}
\times\CP^1$; apart from this, we use the same notations as in \S
\ref{ss:nontoric1}--\ref{ss:nontoric3}.
In particular, it is easy to check that
Proposition \ref{prop:cp2slag} still holds. Hence
we consider the same family of special Lagrangian tori $T_{\gamma(r),
\lambda}$ as above. Lemmas \ref{l:cp2maslov} and \ref{l:cp2wall} also
remain valid, except that the Maslov index formula in Lemma \ref{l:cp2maslov}
becomes \begin{equation}\label{eq:maslovcp1cp1}
\mu(\beta)=2(\beta\cdot[f^{-1}(w)]+\beta\cdot[L^1_\infty]+
\beta\cdot[L^2_\infty]).\end{equation}

In the Clifford case ($r>|\epsilon|$), the superpotential can again be
computed by deforming to the toric case. Denoting again by $z_1$ and $z_2$
the holomorphic coordinates associated to the relative classes $\beta_1$
and $\beta_2$ parallel to the $x$ and $y$ coordinate axes in $(\C^2,
S^1(r_1)\times S^1(r_2))$, we get
\begin{equation}\label{eq:cp1cp1clifford}
W=z_1+z_2+\frac{e^{-\Lambda_1}}{z_1}+\frac{e^{-\Lambda_2}}{z_2},
\end{equation}
where $\Lambda_i$ are the symplectic areas of the two $\CP^1$ factors.
(For simplicity we are only considering the special case $\Lambda_1=\Lambda_2$,
but we keep distinct notations in order to hint at the general case).

On the Chekanov side ($r<|\epsilon|$), we analyze holomorphic discs in
$(\CP^1\times\CP^1,T_{\gamma,0})$ similarly to the case of
$\CP^2$. We denote again by $\beta$ the class of the trivial section of $f$
over the disc $\Delta$ bounded by $\gamma$, and by $\alpha$ the class of
the Lefschetz thimble; and we denote by $H_1=[\CP^1\times\{pt\}]$ and
$H_2=[\{pt\}\times\CP^1]$. Then we have:

\begin{proposition} The only classes in $\pi_2(\CP^1\times\CP^1,T_{\gamma,0})$
which may contain holomorphic discs of Maslov index 2 are $\beta$,
$H_1-\beta-\alpha$, $H_1-\beta$, $H_2-\beta$, and $H_2-\beta+\alpha$.
Moreover, $T_{\gamma,0}$ bounds a unique $S^1$-family of holomorphic discs
in each of these classes, and the corresponding evaluation maps all have
degree 1.
\end{proposition}

\proof We compute the intersection numbers of $\alpha$, $\beta$, $H_1$ and
$H_2$ with the coordinate axes, the fiber $f^{-1}(\epsilon)$,
and the lines at infinity:
\begin{center}
\begin{tabular}{|c||c|c|c|c|c||c|}
\hline
class&$x$-axis&$y$-axis&$L^1_\infty$&$L^2_\infty$&$f^{-1}(\epsilon)\vphantom{\Big|}$&\ $\mu$\ \ \\
\hline
$\alpha\vphantom{\Big|}$&\!\!$-1$&1&0&0&0&0\\
$\beta$&\,\,0&0&0&0&1&2\\
$H_1\vphantom{\Big|}$&\,\,0&1&0&1&1&4\\
$H_{2}\vphantom{{}_\big|}$&\,\,1&0&1&0&1&4\\
\hline
\end{tabular}
\end{center}
The Maslov index formula (\ref{eq:maslovcp1cp1}) and positivity of
intersections with $f^{-1}(\epsilon)$, $L^1_\infty$ and $L^2_\infty$
imply that a holomorphic disc of Maslov index 2 must represent one of
the classes $\beta+k\alpha$, $H_1-\beta+k\alpha$, or $H_2-\beta+k\alpha$,
for some $k\in\Z$. Positivity of intersections with the $x$ and $y$ axes
further restricts the list to the five possibilities mentioned in the
statement of the proposition.

Discs in the class $\beta$ are sections of $f$ over the disc $\Delta$
bounded by $\gamma$; since they are contained in $\C^2$, they are the same
as in the case of $\CP^2$. Discs in the other classes are sections of $f$
over the complement $\CP^1\setminus\Delta$. Denote
by $\varphi$ the biholomorphism from $D^2$ to $\CP^1\setminus\Delta$
such that $\varphi(0)=\infty$ and $\varphi(a)=0$ for some $a\in(0,1)$:
we are looking for holomorphic maps $u:D^2\to\CP^1\times\CP^1$ such that
$f\circ u=\varphi$ and $u(\partial D^2)\subset T_{\gamma,0}$. Considering
the map $q:(x,y)\mapsto x/y$, we see that $q\circ u$ has either a pole or a zero at
$0$ and at $a$, depending on the class represented by $u$, and takes
non-zero complex values everywhere else; moreover it maps the unit circle
to itself. It follows that $q\circ u$ has degree 2 and can be expressed
as $q\circ u(z)=e^{2i\theta}\,z^{\pm 1}\tau_a(z)^{\pm 1}$, where
$e^{i\theta}\in S^1$ and $\tau_a(z)=(z-a)/(1-az)$.
Choosing a square root $\sqrt{h}$ of $h(z)=z\,\varphi(z)/\tau_a(z)$, we
conclude that $u$ is one of
\begin{eqnarray*}
z\mapsto (e^{i\theta}\,z^{-1}\,\tau_a(z)\,\sqrt{h(z)},\,e^{-i\theta}\sqrt{h(z)}),
&& z\mapsto (e^{i\theta}\,\tau_a(z)\,\sqrt{h(z)},\,e^{-i\theta}\,z^{-1}\sqrt{h(z)}),\\
z\mapsto (e^{i\theta}\,z^{-1}\sqrt{h(z)},\,e^{-i\theta}\,\tau_a(z)\,\sqrt{h(z)}),
&& z\mapsto (e^{i\theta}\sqrt{h(z)},\,e^{-i\theta}\,z^{-1}\,\tau_a(z)\sqrt{h(z)}).
\end{eqnarray*}
\endproof

\noindent As before, this implies:

\begin{corollary}
For $r<|\epsilon|$, the superpotential is given by
\begin{equation}
\label{eq:cp1cp1chekanov}
W=u+\frac{e^{-\Lambda_1}(1+w)}{uw}+\frac{e^{-\Lambda_2}(1+w)}{u},
\end{equation}
where $u$ and $w$ are the coordinates associated to the
classes $\beta$ and $\alpha$ respectively.
\end{corollary}

Comparing the formulas (\ref{eq:cp1cp1clifford}) and
(\ref{eq:cp1cp1chekanov}), we see that they are related by the change of
variables \begin{equation}\label{eq:chvar}u=z_1+z_2,\quad
w=z_1/z_2.\end{equation} As in the case of $\CP^2$, this can
be understood in terms of wall-crossing and quantum corrections; the
discussion almost identical to that in \S \ref{ss:nontoric3} and we omit it.

However, we would like to point out one a slightly disconcerting feature
of this example. Since we have assumed that $\Lambda_1=\Lambda_2=\Lambda$,
the right-hand side of (\ref{eq:cp1cp1chekanov}) simplifies
to $u+e^{-\Lambda}(1+w)^2/uw$; this Laurent polynomial has only two critical
points, instead of four for the right-hand side of
(\ref{eq:cp1cp1clifford}) ($z_1=\pm e^{-\Lambda/2}$, $z_2=\pm
e^{-\Lambda/2}$). In particular, the critical value $0$ is lost in the
change of variables, which is unexpected considering the discussion after
Proposition \ref{prop:wallcross}. The reason is of course that the change
of variables (\ref{eq:chvar}) does not quite map $(\C^*)^2$ to itself, and
the critical points where $z_1+z_2=0$ are missing in the $(u,w)$ picture.

\section{Critical values and quantum cohomology} \label{s:m0c1}
The goal of this section is to discuss a folklore result which asserts
that the critical values of the mirror superpotential are the eigenvalues of
quantum multiplication by $c_1(X)$. The argument we present is known to
various experts in the field (Kontsevich, Seidel, ...), but
to our knowledge it has not appeared in the literature.
We state the result specifically in the toric case; however, 
there is a more general relation between the superpotential and $c_1(X)$,
see Proposition \ref{prop:c1cap} below.

\begin{theorem}\label{thm:c1m0}
Let $X$ be a smooth toric Fano variety, and let $W:M\to\C$ be the mirror
Landau-Ginzburg model. Then all the critical values of\/ $W$ are eigenvalues
of the linear map $QH^*(X)\to QH^*(X)$ given by quantum
cup-product with $c_1(X)$.
\end{theorem}

\subsection{Quantum cap action on Lagrangian Floer homology}
The key ingredient in the proof of Theorem \ref{thm:c1m0} is the {\it
quantum cap action} of the quantum cohomology of $X$ on Lagrangian Floer homology. While
the idea of quantum cap action on Floer homology of symplectomorphisms
essentially goes back to Floer \cite{Fl}, its counterpart in the
Lagrangian setting has been much less studied; it can be viewed as a
special case of
Seidel's construction of open-closed operations on Lagrangian Floer homology
(see e.g.\ Section 4 of \cite{SeVCM2}). We review the construction, following
ideas of Seidel and focusing on the specific setting that
will be relevant for Theorem \ref{thm:c1m0}.

The reader is also referred to Biran and Cornea's work \cite{BC}, which
gives a very detailed and careful account of this construction using a slightly
different approach.

Let $L$ be a compact oriented, relatively spin Lagrangian submanifold in
a compact symplectic manifold $(X^{2n},\omega)$ equipped with an almost-complex
structure $J$, and let $\nabla$ be a flat $U(1)$-connection on the trivial
bundle over $L$. We start by describing the operation at the chain level.
Following Fukaya-Oh-Ohta-Ono \cite{FO3}, we use singular chains as the
starting point for the Floer complex, except we use complex coefficients and
assume convergence of all power series. Moreover, for simplicity we quotient
out by those chains whose support is contained in that of a lower-dimensional
chain (this amounts to treating pseudocycles as honest cycles, and allows us
to discard irrelevant terms even when working at the chain level).

Given a class $\beta\in\pi_2(X,L)$, we denote by
$\hat{\mathcal{M}}(L,\beta)$ the space of $J$-holomorphic maps from
$(D^2,\partial D^2)$ to $(X,L)$ representing the class $\beta$
(without quotienting by automorphisms of the disc).
We denote by 
$\hat{ev}_{\beta,\pm 1}:\hat{\mathcal{M}}(L,\beta)\to L$ and 
$\hat{ev}_{\beta,0}:\hat{\mathcal{M}}(L,\beta)\to X$ the evaluation maps at
the boundary points $\pm 1\in \partial D^2$ and the interior point $0\in D^2$.
(So in fact we think of $\hat{\mathcal{M}}(L,\beta)$ as a moduli space
of pseudoholomorphic discs with two marked points on the boundary and one
marked point in the interior, constrained to lie on the geodesic between the
two boundary marked points). 

We will assume throughout this section that
the spaces $\hat{\mathcal{M}}(L,\beta)$ carry well-defined fundamental
chains (of dimension $n+\mu(\beta)$),
and that the evaluation maps are transverse to the chains in $L$ and $X$
that we consider; typically it is necessary to introduce suitable
perturbations in order for these assumptions to hold, but none will be
needed for the application that we have in mind.

\begin{definition}\label{def:qcap}
Let $C\in C_*(L)$ and $Q\in C_*(X)$ be chains in $L$ and $X$ respectively,
such that $C\times Q$ is transverse to the evaluation maps $\hat{ev}_{\beta,1}\times
\hat{ev}_{\beta,0}:\hat{\mathcal{M}}(L,\beta)\to L\times X$. Then we define
\begin{equation}\label{eq:qcap}
Q\cap C=\!\sum_{\beta\in \pi_2(X,L)}z_\beta\ Q\cap_\beta C\in C_*(L),
\end{equation}
where $z_\beta=\exp(-{\textstyle\int_\beta\omega})\hol_\nabla(\partial\beta)$
and $$Q\cap_\beta C=(\hat{ev}_{\beta,-1})_*(\hat{ev}_{\beta,1}\times
\hat{ev}_{\beta,0})^*(C\times Q).$$
\end{definition}

\noindent
In terms of the cohomological degrees $\deg(C)=n-\dim C$ and $\deg(Q)=2n-
\dim Q$, the term $Q\cap_\beta C$
has degree $\deg(C)+\deg(Q)-\mu(\beta)$.

Recall that the Floer differential $\delta=\m_1:C_*(L)\to C_*(L)$
is defined in terms of the moduli spaces
of pseudoholomorphic discs with two marked points on the boundary,
$\mathcal{M}_2(L,\beta)=\hat{\mathcal{M}}(L,\beta)/\R$ (where $\R$
is the stabilizer of $\{\pm 1\}$), and the corresponding evaluation
maps $ev_{\beta,\pm 1}:\mathcal{M}_2(L,\beta)\to L$, by the formula
$$\delta(C)=\partial C+\sum_{\beta\neq 0}z_\beta\,\delta_\beta(C),\ \ \mathrm{where}
\ \ \delta_\beta(C)=(ev_{\beta,-1})_*(ev_{\beta,1})^*(C).$$
We denote by $\delta'(C)=\sum_{\beta\neq 0}z_\beta\delta_\beta(C)$ the
``quantum'' part of the differential.

Assuming there is no obstruction to the construction of Floer
homology, the cap product (\ref{eq:qcap}) descends
to a well-defined map $$\cap:H_*(X)\otimes HF(L,L)\to HF(L,L).$$
In general, the failure of the cap product to be a chain map is encoded
by a higher order operation defined as follows. Let 
$\hat{\mathcal{M}}_3^+(L,\beta) \simeq \hat{\mathcal{M}}(L,\beta)\times\R$
be the moduli space of $J$-holomorphic maps from $(D^2,\partial D^2)$ to
$(X,L)$, with one interior marked point at $0$ and three marked points on
the boundary at $\pm 1$ and at $q=\exp(i\theta)$, $\theta\in (0,\pi)$.
We denote by $\hat{ev}^+_{\beta,q}:\hat{\mathcal{M}}_3^+(L,\beta)\to L$ the
evaluation map at the extra marked point. Define similarly the moduli
space $\hat{\mathcal{M}}_3^-(L,\beta)$ of pseudoholomorpic discs with an
extra marked point at $q=\exp(i\theta)$, $\theta\in(-\pi,0)$, and the
evaluation map $\hat{ev}^-_{\beta,q}$. Then given chains $C,C'\in
C_*(L)$ and $Q\in C_*(X)$ in transverse position, we define
$$\mathfrak{h}^\pm(C,C',Q)=\sum_{\beta\in\pi_2(X,L)} z_\beta\ 
\mathfrak{h}^\pm_\beta(C,C',Q),$$ 
where \hfil $\mathfrak{h}^\pm_\beta(C,C',Q)=(\hat{ev}_{\beta,-1})_*(\hat{ev}_{\beta,1}
\times\hat{ev}^\pm_{\beta,q}\times \hat{ev}_{\beta,0})^*(C\times
C'\times Q).$\hfil\hfil\medskip

Note that the term $\mathfrak{h}^\pm_\beta(C,C',Q)$ has degree
$\deg(C)+\deg(C')+\deg(Q)-\mu(\beta)-1$.
Also recall from \S \ref{s:W} that the obstruction $\m_0\in C_*(L)$ is defined
by $$\m_0=\sum_{\beta\neq 0} z_\beta\
(ev_{\beta,1})_*[\mathcal{M}_1(L,\beta)],$$
where $\mathcal{M}_1(L,\beta)$ is the moduli space of pseudoholomorphic
discs in the class $\beta$ with a single boundary marked point.

\begin{proposition}\label{prop:qcapchain}
Assume that all the chains are transverse to the appropriate
evaluation maps. Then up to signs we have
\begin{equation}\label{eq:qcapchain}
\delta(Q\cap C)=\pm(\partial Q)\cap C \pm Q\cap \delta(C) \pm
\mathfrak{h}^+(C,\m_0,Q)\pm \mathfrak{h}^-(C,\m_0,Q).\end{equation}
\end{proposition}

\proof[Sketch of proof]
The boundary $\partial(Q\cap C)$ of the chain $Q\cap C$ consists of several
pieces, corresponding to the various possible limit scenarios:
\begin{enumerate}
\item One of the two input marked points is mapped to the boundary of the
chain on which it is constrained to lie. The corresponding terms are
$(\partial Q)\cap C$ and $Q\cap (\partial C)$.\smallskip
\item Bubbling occurs at one of the boundary marked points (equivalently, after
reparametrizing this corresponds to the situation where the
interior marked point which maps to $Q$ converges to one of the
boundary marked points). The case where the bubbling occurs at the
incoming marked point $+1$ yields a term $Q\cap \delta'(C)$, while
the case where the bubbling occurs at the outgoing marked point $-1$
yields a term $\delta'(Q\cap C)$. \smallskip
\item Bubbling occurs at some other point of the boundary of the disc, 
i.e.\ we reach the boundary of $\hat{\mathcal{M}}(L,\beta)$; the resulting
contributions are $\mathfrak{h}^+(C,\m_0,Q)$ when bubbling occurs along the
upper half of the boundary, and $\mathfrak{h}^-(C,\m_0,Q)$ when it occurs
along the lower half.
\end{enumerate}
\vskip-1em
\endproof

We will consider specifically the case where $L$ does not bound any 
non-constant pseudoholomorphic discs of Maslov index less than 2; then
the following two lemmas show that Floer homology and the quantum cap action
are well-defined. (In fact, it is clear from the arguments that the relevant
property is the fact that $\m_0$ is a scalar multiple of the fundamental
class $[L]$). The following statement is part of the machinery
developed by Fukaya, Oh, Ohta and Ono \cite{FO3} (see also \cite{Cho1}):

\begin{lemma}\label{l:floerwelldef}
Assume that $L$ does not bound any non-constant pseudoholomorphic discs of
Maslov index less than $2$. Then:
\begin{enumerate}
\item $\m_0$ is a scalar multiple of the fundamental class $[L]$;
\item the Floer cohomology $HF(L,L)$ is well-defined, and $[L]$ is a Floer
cocycle;
\item the chain-level product $\m_2$ determines a well-defined associative
product on $HF(L,L)$, for which $[L]$ is a unit.
\end{enumerate}
\end{lemma}

\proof[Sketch of proof]
(1) The virtual dimension of $\mathcal{M}_1(L,\beta)$ is $n-2+\mu(\beta)$,
so for degree reasons the only non-trivial contributions to $\m_0$ come from
classes of Maslov index 2, and $\m_0$ is an $n$-chain; moreover,
minimality of the Maslov index precludes bubbling, so that $\m_0$ is
actually a cycle, i.e.\ a scalar multiple of the fundamental class $[L]$.

(2) It is a well-known fact in Floer theory (see e.g.\ \cite{FO3}) that
the operations $(\m_k)_{k\ge 0}$ satisfy the $A_\infty$ equations. In
particular, for all $C\in C_*(L)$ we have \begin{equation}\label{eq:m1sq}
\m_1(\m_1(C))+\m_2(C,\m_0)+(-1)^{\deg(C)+1}\m_2(\m_0,C)=0.\end{equation}
To prove that $\m_1\,(=\delta)$ squares to zero,
it is enough to show that \begin{equation}\label{eq:m2unit}
\m_2(C,[L])=C\ \ \text{and}\ \ \m_2([L],C)=(-1)^{\deg(C)}\,C,\end{equation}
since it then follows that the last two terms in (\ref{eq:m1sq})
cancel each other. 

Recall that the products $\m_2(C,[L])$ and $\m_2([L],C)$ are
defined by considering for each class $\beta\in \pi_2(X,L)$ the moduli
space of $J$-holomorphic discs with three boundary marked points in the
class $\beta$, requiring two of the marked points to map to $C$ and $[L]$
respectively, and evaluating at the third marked point. However, the
incidence condition corresponding to the chain $[L]$ is vacuous; so,
provided that $\beta\neq 0$, by forgetting
the unnecessary marked point we see that the construction yields
a chain whose support is contained in that of $\delta_\beta(C)$, which has dimension one
less. It follows that the only nontrivial contribution comes from constant
discs; that contribution is precisely $C$, up to a sign factor left to the reader.

The fact that $[L]$ is a Floer cocycle follows from the observation that,
for any relative class $\beta\neq 0$ containing holomorphic discs, $\delta_\beta([L])$ is a 
chain of dimension $n-1+\mu(\beta)\ge n+1$ in $L$, and hence trivial.
It follows that $\delta([L])=\partial[L]=0$.

(3) In the unobstructed case, the compatibility of $\m_2$ with the Floer differential and its
associativity at the level of cohomology follow from the
$A_\infty$ equations. When $\m_0$ is nonzero, terms of the form
$\m_k(\dots,\m_0,\dots)$ appear
in these equations and make them harder to interpret geometrically.
However, we observe that $\m_k(\dots,[L],\dots)\equiv 0$ for all $k\ge 3$.
Indeed, $\m_k(a_{k-1},\dots,a_i,[L],a_{i-1},\dots,a_1)$ counts
pseudoholomorphic discs with $k+1$ boundary
marked points, where the incidence condition at one of the marked points
is given by the fundamental cycle $[L]$ and hence vacuous; as above,
deleting this extraneous marked point shows that the contribution of each
relative class $\beta\in\pi_2(X,L)$ to
$\m_k(a_{k-1},\dots,a_i,[L],a_{i-1},\dots,a_1)$
has support contained in that of the corresponding contribution to
$\m_{k-1}(a_{k-1},\dots,a_1)$, which has dimension one less.
Hence we can ignore all the terms involving $\m_0$, and the properties
of $\m_2$ are the same as in the unobstructed case.

Finally, the fact that $[L]$ is a unit for the product in Floer cohomology
follows directly from (\ref{eq:m2unit}) (recalling that the sign conventions
in $A_\infty$-algebras are different from those of usual differential graded
algebras, see e.g.\ \cite{FO3,sebook}).
\endproof

\begin{lemma}\label{l:qcapwelldef}
Assume that $L$ does not bound any non-constant pseudoholomorphic discs of
Maslov index less than $2$.
Then the cap product descends to a well-defined
map $\cap:H_*(X)\otimes HF(L,L)\to HF(L,L)$.
\end{lemma}

\proof
By the previous lemma, $\m_0$ is a scalar multiple of $[L]$. Next we observe
that, in the construction of $\mathfrak{h}^\pm(C,[L],Q)$, the
incidence constraint at the extra marked point $q=\exp(i\theta)$ is vacuous.
So the support of each chain $\mathfrak{h}^\pm_\beta(C,[L],Q)$ is
contained in that of the chain $Q\cap_\beta C$, which has dimension one less.
This allows us to discard the terms $\mathfrak{h}^\pm(C,\m_0,Q)$ in 
(\ref{eq:qcapchain}).

Therefore, Proposition \ref{prop:qcapchain} implies that, if $Q$ is a
cycle ($\partial Q=0$) and $C$ is a Floer cocycle ($\delta(C)=0$), then
$Q\cap C$ determines a Floer cocycle, whose class depends only on the
classes of $Q$ and $C$.
\endproof

Next we show that the cap product makes $HF(L,L)$ a module over the
quantum cohomology ring of $X$. We denote by $*$ the quantum cup-product
on $QH^*(X)=H^*(X,\C)$, working again with complex coefficients, i.e.,
specializing the Novikov parameters appropriately so that $J$-holomorphic
spheres in a class $A\in H_2(X)$ are counted with a coefficient 
$\exp(-\int_A\omega)$, and assuming convergence as usual. Moreover, we
use Poincar\'e duality and work with homology instead of cohomology.

\begin{proposition}\label{prop:qhmodule}
Assume that $\m_0$ is a multiple of $[L]$, so that
the cap product $\cap:H_*(X)\otimes HF(L,L)\to HF(L,L)$ is 
well-defined. Then for any $[C]\in HF(L,L)$, $$[X]\cap[C]=[C],$$
and for any $[Q_1],[Q_2]\in H_*(X)$, $$[Q_1]\cap([Q_2]\cap
[C])=([Q_1]*[Q_2])\cap[C].$$
\end{proposition}

\proof[Sketch of proof] We first show that $[X]$ acts by identity. Observe
that, in the definition of the cap product, the incidence constraint at the interior
marked point $0\in D^2$ is vacuous when $Q=[X]$. So for $\beta\neq 0$ the support of
the chain $[X]\cap_\beta C$ is contained in that of the chain
$\delta_\beta(C)$, which has dimension one less. Hence, nonconstant
holomorphic discs contribute trivially to $[X]\cap [C]$. On the other hand
the contribution of constant discs is just the classical intersection of
chains, so that $[X]\cap [C]=[C]$.

To prove the second part of the proposition, consider the moduli space
$\hat{\mathcal{M}}^{(2)}(L,\beta)$ of $J$-holomorphic maps from
$(D^2,\partial D^2)$ to $(X,L)$ with two boundary marked points
at $\pm 1$ and two interior marked points on the real axis, at
$-1<q_1<q_2<1$ (up to simultaneous translation). Denote by
$\hat{ev}_{\beta,q_1}$ and $\hat{ev}_{\beta,q_2}$ the evaluation maps
at the interior marked points, and define
$$\Theta(Q_1,Q_2,C)=\sum_{\beta\in\pi_2(X,L)} z_\beta\
(\hat{ev}_{\beta,-1})_*(\hat{ev}_{\beta,1}\times \hat{ev}_{\beta,q_1}
\times \hat{ev}_{\beta,q_2})^*(C\times Q_1\times Q_2).$$
Given representatives $Q_1, Q_2$ of the given classes $[Q_1],[Q_2]$
and a chain $C\in C_*(L)$, assuming
transversality as usual, a case-by-case
analysis similar to the proof of Proposition \ref{prop:qcapchain} shows that
$$\delta(\Theta(Q_1,Q_2,C))=\pm \Theta(Q_1,Q_2,\delta(C))\pm Q_1\cap
(Q_2\cap C)\pm (Q_1*Q_2)\cap C\pm (\m_0\mbox{-terms}).$$
More precisely, the boundary of the chain $\Theta(Q_1,Q_2,C)$ consists of:
\begin{enumerate}
\item $\Theta(Q_1,Q_2,\partial C)$, corresponding to the situation where
the input marked point at $+1$ is mapped to the boundary of the chain
$C$. (Note that since $Q_1$ and $Q_2$ are cycles, we do not include the two
terms $\Theta(\partial Q_1,Q_2,C)$ and $\Theta(Q_1,\partial Q_2,C)$ which
would be present in the general case.)\smallskip
\item $\delta'(\Theta(Q_1,Q_2,C))$ and $\Theta(Q_1,Q_2,\delta'(C))$,
corresponding to bubbling at one of the boundary marked points $\pm 1$
(or equivalently after reparametrization, the situation where the interior
marked points $q_1,q_2$ both converge to $\pm 1$).
\smallskip
\item $Q_1\cap (Q_2\cap C)$, corresponding to the situation where $q_1\to
-1$ (or equivalently up to reparametrization, $q_2\to 1$), resulting in a
two-component map with one interior marked point in each component.
\smallskip
\item $(Q_1*Q_2)\cap C$, corresponding to the situation where the two marked
points $q_1$ and $q_2$ come together, leading to the bubbling of a sphere
component which carries both marked points, attached to
the disc at a point on the real axis.
\smallskip
\item terms involving $\m_0$ and higher order operations defined
analogously to $\mathfrak{h}^\pm$, involving
moduli spaces of discs with three marked points on the boundary and two
in the interior; these occur as in the proof of Proposition
\ref{prop:qcapchain} when bubbling occurs at a point of $\partial
D^2\setminus\{\pm 1\}$. 
\end{enumerate}
By the same argument as in the proof of Lemma \ref{l:qcapwelldef}, when
$\m_0$ is a multiple of $[L]$ we can safely ignore the last set of terms
because the corresponding chains are supported on lower-dimensional subsets.
Thus, if we assume that $C$ is a Floer cocycle (i.e., $\delta(C)=0$), the
above formula implies that the Floer cocycles $Q_1\cap(Q_2\cap C)$ and
$(Q_1*Q_2)\cap C$ represent the same Floer cohomology class (up to signs,
which are left to the reader).
\endproof

\subsection{Cap product by $c_1(X)$ and proof of Theorem
\ref{thm:c1m0}}

Let $(X,\omega,J)$ be a smooth compact K\"ahler manifold of complex
dimension $n$, equipped with a holomorphic $n$-form $\Omega$ defined over
the complement of an anticanonical divisor $D$. 
Let $L\subset X\setminus D$ be a special
Lagrangian submanifold, or more generally a Lagrangian submanifold whose
Maslov class vanishes in $X\setminus D$ (so that Lemma \ref{l:maslov}
holds), and let $\nabla$ be a flat $U(1)$-connection on the trivial line
bundle over $L$.
We assume that $L$ does not bound any nonconstant holomorphic disc
of Maslov index less than 2, so that the Floer obstruction $\m_0$ is a
constant multiple of the fundamental class, $\m_0=m_0(L,\nabla)\,[L]$,
and Floer homology and the quantum cap product are well-defined.

\begin{lemma}\label{l:c1cap}
$c_1(X)\cap [L]=m_0(L,\nabla)\,[L]$.
\end{lemma}

\proof 
We actually compute the cap product $[D]\cap [L]$. Since $[D]\cap_\beta [L]$
is a chain of dimension $n-2+\mu(\beta)$ in $L$, the only contributions
to $[D]\cap [L]$ come from classes of Maslov index at most $2$. Moreover,
since $L\subset X\setminus D$, there are no contributions from constant
discs, so we only need to consider classes with $\mu(\beta)=2$.

By Lemma \ref{l:maslov} every holomorphic map $u:(D^2,\partial D^2)\to
(X,L)$ of Maslov index 2 intersects $D$ in a single point $u(z_0)$, $z_0\in
D^2$. Moreover, for every $q=e^{i\theta}\in\partial D^2$ there exists a
unique automorphism of $D^2$ which maps $z_0$ to $0$ and $q$ to $-1$.
It follows that $[D]\cap_\beta [L]$ is the chain consisting of all
boundary points of all holomorphic discs in the class $\beta$, i.e.\ 
$[D]\cap_\beta [L]=(ev_{\beta,1})_*[\mathcal{M}_1(L,\beta)]$. Summing
over $\beta$, we conclude that $[D]\cap [L]=\m_0=m_0(L,\nabla)\,[L]$.
\endproof

Lemma \ref{l:c1cap} implies the following proposition, which is the
core of Theorem \ref{thm:c1m0}.

\begin{proposition}\label{prop:c1cap}
If $HF(L,L)\neq 0$ then $m_0(L,\nabla)$ is an eigenvalue of the linear
map $\Lambda:QH^*(X)\to QH^*(X)$ defined by $\Lambda(\alpha)=\alpha*c_1(X)$.
\end{proposition}

\proof Since $[L]$ is the unit for the product on Floer cohomology,
the assumption $HF(L,L)\neq 0$ implies that $[L]$ is a nonzero element
of $HF(L,L)$. Lemma \ref{l:c1cap} states that $(c_1(X)-m_0(L,\nabla))
\cap [L]=0$. But then Proposition \ref{prop:qhmodule} implies that
quantum cup-product by $c_1(X)-m_0(L,\nabla)$ is not invertible. \endproof

The only remaining ingredient of the proof of Theorem \ref{thm:c1m0} is to
show that critical points of the superpotential correspond to special
Lagrangians with nonzero Floer homology. This follows from a general
interpretation of the $k$-th derivatives of the superpotential in terms
of $\m_k$, at least in the toric Fano case (see the works of Cho and Oh
\cite{ChoOh, Cho}). Theorem 10.1 in \cite{ChoOh} states that, if $X$ is
a toric Fano variety and $L$ is a toric fiber, then $HF(L,L)$ is non-zero
if and only if the contributions of Maslov index 2 classes to $\delta([pt])$
cancel out in $H_1(L)$. In our terminology, the statement is:

\begin{proposition}[Cho-Oh \cite{ChoOh}]\label{prop:hfm1}
Let $L$ be a toric fiber in a toric Fano variety, equipped with a flat
$U(1)$ connection $\nabla$. Then $HF(L,L)\neq 0$ if and only if
\begin{equation}\label{eq:m1zero}
m_1(L,\nabla):=\sum_{\mu(\beta)=2} n_\beta(L) \exp(-{\textstyle \int_\beta\omega})
\hol_\nabla(\partial\beta)\,[\partial\beta]=0\in H_1(L,\C).\end{equation}
\end{proposition}

The ``only if'' part actually holds even in the non-toric case, assuming
the minimal Maslov index of a holomorphic disc to be 2. Indeed,
it is easy to check that
if $C$ is a codimension 1 cycle in $L$ (i.e.\ $\partial C=0$) then \begin{equation}
\label{eq:deltacodim1}\delta(C)=\pm ([C]\cdot m_1(L,\nabla))\,[L]\end{equation}
(classes of Maslov index $>2$ do not contribute to $\delta(C)$ for dimension
reasons), so that when $m_1(L,\nabla)\neq 0$ the fundamental class $[L]$ can be
realized as a Floer coboundary. However, to our knowledge the ``if'' part
of the statement has not been proved outside of the toric case; the argument in
\cite{ChoOh} relies on the specific features of holomorphic discs in toric
varieties to show that classes of Maslov index $>2$ never contribute to the
Floer differential (Proposition 7.2 in \cite{ChoOh}), so that (\ref{eq:deltacodim1}) holds for cycles of any
dimension, and vanishing of $m_1(L,\nabla)$ implies nontriviality of the Floer
homology.

Finally, recall from Section \ref{s:defs} that $T_{(L,\nabla)}M\simeq 
\mathcal{H}^1_\psi(L)\otimes\C\simeq H^1(L,\C)$, by mapping $(v,\alpha)\in
T_{(L,\nabla)}M\subset C^\infty(NL)\oplus \Omega^1(L,\R)$ to
$[-\iota_v\omega+i\alpha]$. Then we have:

\begin{lemma}\label{l:dW}
The differential of $W=m_0:M\to\C$ is
$$dW_{(L,\nabla)}(v,\alpha)=\langle [-\iota_v\omega+i\alpha],m_1(L,\nabla)
\rangle.$$
\end{lemma}

\proof Let $z_\beta=\exp(-\int_\beta\omega)\hol_\nabla(\partial\beta)$,
and observe as in the proof of Lemma \ref{l:holomcoords} that 
$d\log z_\beta(v,\alpha)=\langle [-\iota_v\omega+i\alpha],
[\partial\beta]\rangle$ (by Stokes' theorem). Hence, the differential of
$W=\sum n_\beta(L)\,z_\beta$ is $dW(v,\alpha)=\sum n_\beta(L)\,z_\beta\,
\langle [-\iota_v\omega+i\alpha],[\partial\beta]\rangle$.
\endproof

Theorem \ref{thm:c1m0} now follows from Proposition \ref{prop:c1cap},
Proposition \ref{prop:hfm1} and Lemma \ref{l:dW}: 
if $(L,\nabla)$ is a critical point
of $W$ then by Lemma \ref{l:dW} it satisfies $m_1(L,\nabla)=0$, and hence
by Proposition \ref{prop:hfm1} the Floer cohomology $HF(L,L)$ is nontrivial.
Proposition \ref{prop:c1cap} then implies that the critical value
$W(L,\nabla)=m_0(L,\nabla)$ is an eigenvalue of quantum multiplication by
$c_1(X)$.

\section{Admissible Lagrangians and the reference fiber} \label{s:fiber}

In this section we give a brief, conjectural discussion of the manner in
which the mirror construction discussed in the preceding sections relates
to mirror symmetry for the Calabi-Yau hypersurface $D\subset X$. For
simplicity, unless otherwise specified we assume throughout this section
that $D$ is smooth. 

\subsection{The boundary of $M$ and the reference fiber}
Denote by $\sigma\in H^0(X,K_X^{-1})$ the defining section of $D$, and 
identify a tubular neighborhood
$U$ of $D$ with a neighborhood of the zero section in the normal bundle
$N_D\simeq {(K_X^{-1})}_{|D}$, denoting by $p:U\to D$ the projection.
Then we have:

\begin{lemma}\label{l:residue}
$D$ carries a nonvanishing holomorphic $(n-1)$-form $\Omega_D$, called the
{\em residue} of\/ $\Omega$ along $D$, such that, in a neighborhood of $D$,
\begin{equation}\label{eq:residue}
\Omega=\sigma^{-1}d\sigma\wedge p^*\Omega_D+O(1).\end{equation}
\end{lemma}

\noindent
Note that, even though $\sigma^{-1}d\sigma$ depends on the choice of a
holomorphic connection on $K_X^{-1}$ (one can e.g.\ use the Chern
connection), it only does so by a bounded amount, so this ambiguity
has no incidence on (\ref{eq:residue}). The choice of the
projection $p:U\to D$ does not matter either, for a similar reason.

\proof
Near any given point $q\in D$, choose local holomorphic coordinates 
$(x_1,\dots,x_n)$ on $X$ such that $D$ is the hypersurface $x_1=0$.
Then locally we can write $\Omega=x_1^{-1}h(x_1,x_2,\dots,x_n)\,dx_1\wedge
\dots \wedge dx_n$, for some nonvanishing holomorphic function $h$. 
We set $\Omega_D=h(0,x_2,\dots,x_n)\,dx_2\wedge\dots\wedge dx_n$; in
other terms, $\Omega_D=(x_1\iota_{\partial_{x_1}}\Omega)_{|D}$.

If we change the coordinate system to a different one $(y_1,\dots,y_n)$
for which $D$ is again defined by $y_1=0$, then
$x_1=y_1\,\phi(y_1,\dots,y_n)$ for some nonvanishing holomorphic function
$\phi$, so that $x_1^{-1}\,dx_1=y_1^{-1}\,dy_1+d(\log\phi)$. Therefore,
denoting by $\mathcal{J}$ the Jacobian of the change of variables
$(x_2,\dots,x_n)\mapsto (y_2,\dots,y_n)$ on $D$, we have
$x_1^{-1}h\,dx_1\wedge \dots\wedge dx_n=(y_1^{-1}h\,\mathcal{J}+O(1))\,
dy_1\wedge \dots\wedge dy_n$. Hence $\Omega_D$ is well-defined.

Finally, equation (\ref{eq:residue}) follows by considering
a coordinate system in which the first coordinate is exactly $\sigma$ in a
local trivialization of $K_X^{-1}$ and the other coordinates are pulled back
from $D$ by the projection $p$.
\endproof

Lemma \ref{l:residue} shows that $D$ (equipped with the restricted complex
structure and symplectic form, and with the volume form $\Omega_D$) is a
Calabi-Yau manifold.

\begin{remark}\label{rmk:Dnormalcross}
If $D$ has normal crossing singularities, then the same construction
yields a holomorphic $(n-1)$-form $\Omega_D$ which has poles along the
singular locus of $D$. 
\end{remark}

 Assume that $\Lambda$ is a special Lagrangian
submanifold in $(D,\Omega_D)$: then, in favorable cases, we can try to
look for special
Lagrangian submanifolds in $X\setminus D$ which are $S^1$-fibered over
$\Lambda$. For example, if we are in a product situation, i.e.\ locally
$X=\C\times D$, with the product complex structure and K\"ahler form,
and $\Omega=x^{-1}dx\wedge \Omega_D$, then $S^1(r)\times \Lambda$ is
special Lagrangian in $(X\setminus D,\Omega)$.

In general, the classical symplectic neighborhood theorem implies that
$U$ is symplectomorphic to a neighborhood of the zero section in a
symplectic vector bundle over $D$. It is then easy to see that the
preimage of $\Lambda$ is foliated by Lagrangian submanifolds (which
intersect each fiber in concentric circles). Due to the lack of compatibility
between the standard symplectic chart and the holomorphic volume form,
these submanifolds are not special Lagrangian in general, 
but Lemma \ref{l:residue} implies that, as the radius of the circle in
the fiber tends to zero (i.e., as the submanifolds get closer to $D$),
they become closer and closer to being special Lagrangian. 
Thus it is reasonable to hope for the
existence of a nearby family of special Lagrangian submanifolds.

In terms of the moduli space $M$ of pairs $(L,\nabla)$ consisting of a
special Lagrangian submanifold of $(X\setminus D,\Omega)$ and a flat
$U(1)$-connection on the trivial bundle over $L$, this suggests the
following (somewhat optimistic) conjecture:

\begin{conj}\label{conj:boundary}
Near its boundary, $M$ consists of pairs $(L,\nabla)$ such that the
Lagrangian submanifold $L\subset U\cap (X\setminus D)$ is 
a circle bundle over a special Lagrangian submanifold of $D$, with the
additional property that
every fiber bounds a holomorphic disc of Maslov index 2 contained in $U$.
\end{conj}

The main evidence for this conjecture is given by the examples in Sections
\ref{s:toric} and~\ref{s:nontoric} above. In those examples $D$ has normal
crossing singularities, along which $\Omega_D$ has poles, so special
Lagrangian submanifolds of $D$ are defined as in Section \ref{s:defs}.
In this setting, the moduli space $M$ also has corners, corresponding to
the situation where $L$ lies close to the singular locus of $D$.
Apart from this straightforward adaptation, the boundary structure of the
moduli space of special Lagrangians is exactly as described by the conjecture.

More precisely, as one approaches the boundary of $M$, the special Lagrangian submanifold
$L$ collapses to a special Lagrangian submanifold $\Lambda$ of $D$, and the
collapsing makes $L$ a (topologically trivial) $S^1$-bundle over $\Lambda$.
Moreover, each circle fiber bounds a small holomorphic disc which intersects
$D$ transversely in a single point; we denote by $\delta\in \pi_2(X,L)$ the
homotopy class of these discs. As $L$ collapses onto $\Lambda$, the
symplectic area $\int_\delta\omega$ shrinks to zero; in terms of the variable
$z_\delta=\exp(-\int_\delta\omega) \hol_{\nabla}(\partial\delta)$,
we get $|z_\delta|\to 1$. In other terms, Conjecture \ref{conj:boundary}
implies that $\partial M$ is defined by the equation $|z_\delta|=1$.

Among the points of $\partial M$, those where $z_\delta=1$ stand out,
because they correspond to the situation where the holonomy of $\nabla$
is trivial along the fiber of the $S^1$-bundle $L\to\Lambda$, i.e.\ 
$\nabla$ is lifted from a connection on the trivial bundle over $\Lambda$.
The set of such points can therefore be identified with a moduli space
$M_D$ of pairs of special Lagrangian submanifolds in $D$ and flat
$U(1)$-connections over them.

\begin{conj}\label{conj:syzD}
The subset $M_D=\{z_\delta=1\}\subset\partial M$ is the
Strominger-Yau-Zaslow mirror of $D$.
\end{conj}

Assuming these conjectures, it is tempting to think of $M_D$ as the 
{\it reference fiber} (or ``fiber at infinity'') of the Landau-Ginzburg
model $W:M\to\C$. Of course, $M_D$ is not actually a fiber of $W$; but
the contributions of other relative homotopy classes to the superpotential
are negligible compared to $z_\delta$, at least in the rescaling limit
suggested by Conjecture \ref{conj:renormalize}. So, near $\partial M$, we
expect to have $W=z_\delta+o(1)$, and the fiber $W^{-1}(1)$ (when it is
well-defined, which is not always the case considering wall-crossing and
boundary phenomena) can essentially be identified with $M_D$.

Moreover, we expect the boundary of $M$ to fiber over $S^1$. Indeed,
we can set the holonomy of $\nabla$ along the fiber of the $S^1$-bundle
$L\to\Lambda$ to equal any given unit complex number, rather than $1$.
Thus we have:

\begin{proposition}
Assuming Conjecture \ref{conj:boundary}, the map $z_\delta:\partial M\to S^1$
is a fibration with fiber $M_D$.
\end{proposition}

A hasty look at the situation might lead one to conclude, incorrectly,
that the fibration
$z_\delta:\partial M\to S^1$ is trivial. In fact, the symplectic monodromy
of this fibration, viewed as an autoequivalence of the Fukaya category of
$M_D$, is expected to be mirror to the autoequivalence of $D^bCoh(D)$
induced by $\mathcal{E}\mapsto (K_X)_{|D}\otimes \mathcal{E}$.

For example, consider the case where $X=\CP^2$ and $D$ is a smooth elliptic
curve obtained by a small generic deformation of the union of the coordinate
lines. While we do not have an explicit description of $M$ in this setting,
one can use the toric model as a starting point, and the example in Section
\ref{s:nontoric} as an indication of how the smoothing of $D$ affects the
geometry of $M$ near the corners of the moment polytope.
Let $L$ be a special Lagrangian torus which lies close to the portion
of $D$ located along the $x$ and $y$ coordinate axes: then $\delta$ is 
the relative homotopy class called $\beta$ in \S \ref{s:nontoric}; using
the same notations as in \S \ref{s:nontoric}, $\partial\beta$ is the fiber
of the $S^1$-bundle $L\to\Lambda$, while $\partial\alpha$ is a section. 
If we move $\Lambda$ by translations all around the
elliptic curve $D$, and look at the corresponding family of special
Lagrangians $L\subset X\setminus D$, a short calculation shows that
the monodromy acts on $H_1(L)$ by \begin{equation}
\label{eq:monodrD}\partial\beta\mapsto \partial\beta,
\quad \partial\alpha\mapsto \partial\alpha+9\,\partial\beta.\end{equation}
Observe that $M_{D,\theta}=\{z_\beta=e^{i\theta}\}\subset \partial M$ is an
$S^1$-bundle over $S^1$, where the base corresponds to the moduli space 
$\mathcal{B}_D$ of special 
Lagrangians $\Lambda\subset D$ (occurring as collapsed limits of special
Lagrangian tori $L\subset X\setminus D$), and the fiber corresponds
to the holonomy of the flat connection over $L$
with the constraint $\hol_\nabla(\partial\beta)=e^{i\theta}$.
For $\theta=0$, a section of the $S^1$-bundle
$M_{D,0}=M_D\to \mathcal{B}_D$
is given by fixing the holonomy along $\partial\alpha$, e.g.\ requiring
$\hol_\nabla(\partial\alpha)=1$; since $\hol_\nabla(\partial\beta)=1$ this
constraint is compatible with the monodromy (\ref{eq:monodrD}). Now,
increase $\theta$ from $0$ to $2\pi$ and deform this section into
each $M_{D,\theta}=\{z_\beta=e^{i\theta}\}$: it follows from 
(\ref{eq:monodrD}) that we get a section of the $S^1$-bundle $M_{D,\theta}
\to \mathcal{B}_D$ along which the holonomy 
$\hol_\nabla(\partial\alpha)$ varies by $9\,\theta$.
So, when we return to $M_D$ as $\theta$ reaches $2\pi$, the homotopy class
of the section has changed by 9 times the fiber of the $S^1$-bundle
$M_D\to\mathcal{B}_D$: i.e., in this example the monodromy of the fibration
$z_\beta:\partial M\to S^1$ is given by
$$\begin{pmatrix}1&0\\9&1\end{pmatrix}.$$

\subsection{Fukaya categories and restriction functors}
We now return to the general case, and discuss briefly the Fukaya category 
of the Landau-Ginzburg model $W:M\to\C$, assuming that Conjectures 
\ref{conj:boundary} and \ref{conj:syzD} hold. The general idea, which
goes back to Kontsevich \cite{KoENS} and Hori-Iqbal-Vafa \cite{HIV}, is
to allow as objects {\em admissible Lagrangian submanifolds}\/ of $M$, i.e.\
potentially non-compact Lagrangian submanifolds which, outside of a compact
subset, are invariant under the gradient flow of $-\mathrm{Re}(W)$. The case
of Lefschetz fibrations (i.e., when the critical points of $W$ are
nondegenerate) has been studied in great detail by Seidel; in this case,
which is by far the best understood, the theory can be
formulated in terms of the vanishing cycles at the critical points 
(see e.g.\ \cite{sebook}).

The formulation which is the most relevant to us is
the one which appears in Abouzaid's work \cite{Ab1,Ab2}: in this version,
one considers Lagrangian submanifolds of $M$ with boundary contained in a
given fiber of the superpotential, and which near the reference fiber 
are mapped by $W$ to an embedded curve $\gamma\subset\C$. In our case,
using the fact that near $\partial M$ the superpotential is
$W=z_\delta+o(1)$, we consider Lagrangian submanifolds with boundary in
$M_D=\{z_\delta=1\}$:

\begin{definition}
A Lagrangian submanifold $L\subset M$ with (possibly empty)
boundary $\partial L\subset M_D$
is {\em admissible with slope 0} if the restriction of $z_\delta$ to $L$
takes real values near the boundary of $L$.
\end{definition}

Similarly, we say that $L$ is admissible with slope $\theta\in (-\frac\pi2,
\frac\pi2)$ if $\partial L
\subset M_D$ and, near $\partial L$, $z_\delta$ takes values in the
half-line $1-e^{i\theta}\R_+$.
The definition of Floer homology for admissible Lagrangians is the usual
one in this context: to determine $HF(L_1,L_2)$, one first
deforms $L_2$ (rel.\ its boundary) to an admissible Lagrangian $L_2^+$ whose slope is
greater than that of $L_1$, and one computes Floer homology for the pair of
Lagrangians $(L_1,L_2^+)$ inside $M$ (ignoring boundary
intersections).

We denote by $\mathcal{F}(M,M_D)$ the Fukaya category constructed in this
manner. Replacing the superpotential by $z_\delta$ in the definition
has two advantages: on one hand it makes admissibility 
a much more geometric condition, and on the other hand it eliminates
difficulties associated with wall-crossing and definition of the
superpotential. In particular, when comparing the B-model on $X$ and the
A-model on $M$ this allows us to entirely eliminate the superpotential
from the discussion. Since in the rescaling limit of Conjecture
\ref{conj:renormalize} the superpotential is expected to be
$W=z_\delta+o(1)$, we conjecture that
$\mathcal{F}(M,M_D)$ is derived equivalent to the physically relevant
category of Lagrangian submanifolds.

Finally, we conclude the discussion by observing that, by construction, the
boundary of an admissible Lagrangian in $M$ is a Lagrangian submanifold
of $M_D$ (possibly empty, and not necessarily connected). We claim that
there is a well-defined {\it restriction functor} $\rho:\mathcal{F}(M,M_D)\to
\mathcal{F}(M_D)$ from the Fukaya category of $M$ to that of $M_D$, which
at the level of objects is simply $(L,\nabla)\mapsto (\partial
L,\nabla_{|\partial L})$. At the level of morphisms, the restriction functor
essentially projects to the part of the Floer complex generated by the
intersection points near the boundary. More precisely, given an intersection
point $p\in \mathrm{int}(L_1)\cap \mathrm{int}(L_2^+)$, $\rho(p)$ is a
linear combination of intersection points in which the coefficient of
$q\in \partial L_1\cap \partial L_2$ counts the number of holomorphic strips
connecting $p$ to $q$ in $(M,L_1\cup L_2^+)$.

This suggests the following conjecture, which can be thought of as
``relative homological mirror symmetry'' for the pair $(X,D)$:

\begin{conj}\label{conj:relhms}
There is a commutative diagram
$$\begin{CD}
D^bCoh(X) @>\text{restr}>> D^bCoh(D) \\
@V{\simeq}VV @VV{\simeq}V \\
D^\pi\mathcal{F}(M,M_D) @>{\rho}>> D^\pi\mathcal{F}(M_D)
\end{CD}
$$
\end{conj}

In this diagram, the horizontal arrows are the restriction functors, and
the vertical arrows are the equivalences predicted by homological mirror
symmetry.

Some evidence for Conjecture \ref{conj:relhms} is provided by the case of
Del Pezzo surfaces \cite{AKO2}. Even though it is not clear that the
construction of the mirror in \cite{AKO2} corresponds to the one discussed
here, it is striking to observe how the various ingredients fit
together in that example. Namely, by comparing the calculations for Del
Pezzo surfaces in \cite{AKO2} with Polishchuk and Zaslow's work on mirror
symmetry for elliptic curves \cite{PZ}, it is readily apparent that:
\begin{itemize}
\item the fiber of the Landau-Ginzburg model 
$W:M\to\C$ is mirror to an elliptic curve $E$ in the anticanonical
linear system $|K_X^{-1}|$;\smallskip
\item  the Fukaya category of the superpotential admits an exceptional
collection consisting of Lefschetz thimbles; under mirror symmetry for
elliptic curves, their boundaries, which
are the vanishing cycles of the critical points of $W$, correspond exactly
to the restrictions to $E$ of
the elements of an exceptional collection for $D^bCoh(X)$;\smallskip
\item  the behavior of the restriction functors on these exceptional
collections is exactly as predicted by Conjecture \ref{conj:relhms}.
\end{itemize}


\begin{thebibliography}{99}
\bibitem{Ab1}
M. Abouzaid,
{\sl Homogeneous coordinate rings and mirror symmetry for toric varieties},
Geom. Topol. {\bf 10} (2006), 1097--1157 (math.SG/0511644).
\bibitem{Ab2}
M. Abouzaid,
{\sl Morse homology, tropical geometry, and homological mirror symmetry for
toric varieties},
preprint (math.SG/0610004).
\bibitem{AKO1}
D. Auroux, L. Katzarkov, D. Orlov,
{\sl Mirror symmetry for weighted projective planes and their noncommutative
deformations},
to appear in Ann. Math. (math.AG/0404281).
\bibitem{AKO2}
D. Auroux, L. Katzarkov, D. Orlov,
{\sl Mirror symmetry for Del Pezzo surfaces: Vanishing cycles and coherent
sheaves},
Inventiones Math. {\bf 166} (2006), 537--582 (math.AG/0506166).
\bibitem{BC}
P. Biran, O. Cornea,
{\sl Quantum structures for Lagrangian submanifolds}, preprint, in preparation.
\bibitem{chekanov}
Y. Chekanov,
{\sl Lagrangian tori in a symplectic vector space and global 
symplectomorphisms},
Math. Z. {\bf 223} (1996), 547--559.
\bibitem{schlenk}
Y. Chekanov, F. Schlenk, in preparation.
\bibitem{Cho1}
C.-H. Cho, {\sl Holomorphic discs, spin structures and the Floer
cohomology of the Clifford torus}, Int. Math. Res. Not. (2004), 
1803--1843 (math.SG/0308224).
\bibitem{Cho}
C.-H. Cho,
{\sl Products of Floer cohomology of torus fibers in toric Fano manifolds},
Comm. Math. Phys. {\bf 260} (2005), 613--640 (math.SG/0412414).
\bibitem{ChoOh}
C.-H. Cho, Y.-G. Oh,
{\sl Floer cohomology and disc instantons of Lagrangian torus fibers in
Fano toric manifolds}, Asian J. Math. {\bf 10} (2006), 773--814
(math.SG/0308225).
\bibitem{CL}
O. Cornea, F. Lalonde,
{\sl Cluster homology}, preprint (math.SG/0508345).
\bibitem{EP}
Y. Eliashberg, L. Polterovich, 
{\sl The problem of Lagrangian knots in four-manifolds},
Geometric Topology (Athens, 1993), AMS/IP Stud. Adv. Math., Amer. Math.
Soc., 1997, pp. 313--327.
\bibitem{Fl}
A. Floer,
{\sl Symplectic fixed points and holomorphic spheres},
Commun.\ Math.\ Phys.\ {\bf 120} (1989), 575--611.
\bibitem{FO3}
K. Fukaya, Y.-G. Oh, H. Ohta, K. Ono,
{\sl Lagrangian intersection Floer theory: Anomaly and obstruction},
preprint, second expanded version, 2006.
\bibitem{GS}
M. Gross, B. Siebert,
{\sl From real affine geometry to complex geometry},
preprint (math.AG/ 0703822).
\bibitem{hitchin}
N. Hitchin,
{\sl The moduli space of special Lagrangian submanifolds},
Ann.\ Scuola Norm.\ Sup.\ Pisa Cl.\ Sci.\ {\bf 25} (1997), 503--515
(dg-ga/9711002).
\bibitem{HIV}
K.\ Hori, A.\ Iqbal, C.\ Vafa, {\sl D-branes and mirror symmetry},
preprint (hep-th/0005247).
\bibitem{HV}
K. Hori, C. Vafa, 
{\sl Mirror symmetry}, preprint (hep-th/0002222).
\bibitem{KoENS}
M.\ Kontsevich, {\sl Lectures at ENS, Paris, Spring 1998}, notes taken by
J.\ Bellaiche, J.-F.\ Dat, I.\ Marin, G.\ Racinet and H.\ Randriambololona.
\bibitem{KS1}
M. Kontsevich, Y. Soibelman,
{\sl Homological mirror symmetry and torus fibrations},
Symplectic geometry and mirror symmetry (Seoul, 2000), 
World Sci. Publ., 2001, pp. 203--263 (math.SG/0011041).
\bibitem{KS2}
M. Kontsevich, Y. Soibelman,
{\sl Affine structures and non-Archimedean analytic spaces},
The unity of mathematics, Progr. Math. {\bf 244}, Birkh\"auser Boston, 
2006, pp. 321--385 (math.AG/ 0406564).
\bibitem{mclean}
R.\,C. McLean,
{\sl Deformations of calibrated submanifolds},
Comm. Anal. Geom. {\bf 6} (1998), 705--747.
\bibitem{PZ}
A. Polishchuk, E. Zaslow,
{\sl Categorical mirror symmetry: the elliptic curve},
Adv.\ Theor.\ Math.\ Phys.\ {\bf 2} (1998), 443--470.
\bibitem{salur}
S.\ Salur, {\sl Deformations of special Lagrangian submanifolds},
Comm.\ Contemp.\ Math.\ {\bf 2} (2000), 365--372 (math.DG/9906048).
\bibitem{SeGraded}
P. Seidel,
{\sl Graded Lagrangian submanifolds},
Bull. Soc. Math. Fr. {\bf 128} (2000), 103--149.
\bibitem{SeVCM2}
P. Seidel, {\sl More about vanishing cycles and mutation},
``Symplectic Geometry and Mirror Symmetry'', Proc. 4th KIAS International
Conference (Seoul, 2000), World Sci., Singapore, 2001, pp. 429--465
(math.SG/0010032).
\bibitem{seidel}
P. Seidel, 
{\sl A biased view of symplectic cohomology},
preprint (math.SG/0704.2055).
\bibitem{sebook}
P. Seidel,
{\sl Fukaya categories and Picard-Lefschetz theory},
in preparation.
\bibitem{SYZ}
A. Strominger, S.-T. Yau, E. Zaslow,
{\sl Mirror symmetry is T-duality},
Nucl. Phys. B {\bf 479} (1996), 243--259 (hep-th/9606040).
\end{thebibliography}
\end{document}